\documentclass[10pt,a4paper]{article}
\usepackage[utf8]{inputenc}
\usepackage[english]{babel}
\usepackage{amsmath,bm}
\usepackage{amsfonts}
\usepackage{amssymb}
\usepackage{amsthm}
\usepackage{graphicx}
\usepackage[left=2cm,right=2cm,top=2cm,bottom=2cm]{geometry}

\usepackage{color}
\usepackage{tabularx}  
\usepackage{placeins} 

\def\v{\mathbf{v}}
\def\u{\mathbf{u}}
\def\U{\mathbf{U}}
\def\x{\mathbf{x}}
\def\f{\mathbf{f}}
\def\n{\mathbf{n}}

\def\Y{\mathbf{Y}}
\def\T{\mathcal{T}}
\def\VDG{\mathbf{V}_{DG}}

\def\sig{\boldsymbol{\sigma}}
\def\tin{\ \mathrm{d}}

\newtheorem*{remark}{Remark}

\begin{document}

\markboth{Esterhazy, Schneider, Mazzieri, Bokelmann}
{Insights into the modeling of seismic waves}

\title{
INSIGHTS INTO THE MODELING OF SEISMIC WAVES FOR THE DETECTION OF UNDERGROUND CAVITIES
}

\author{SofiEsterhazy, Felix Schneider, Ilario Mazzieri, Götz Bokelmann
}

%
%
%
%
%
%

\maketitle

\begin{abstract}
Motivated by the need to detect an underground cavity within the 
procedure of an On-Site-Inspection (OSI), of the Comprehensive  
Nuclear Test Ban Treaty Organization, the aim of this paper is to 
present results on the comparison of our numerical simulations with 
an analytic solution. The accurate numerical modeling can facilitate 
the development of proper analysis techniques to detect the remnants 
of an underground nuclear test. The larger goal is to help set a 
rigorous scientific base of OSI and to contribute to bringing the 
Treaty into force. For our 3D numerical simulations, we use the 
discontinuous Galerkin Spectral Element Code SPEED jointly developed 
at MOX (The Laboratory for Modeling and Scientific Computing, 
Department of Mathematics) and at DICA (Department of Civil and 
Environmental Engineering) of  the Politecnico di Milano. 
\end{abstract}

\section{Introduction}
If a suspicious seismic signal has been recorded by the International Monitoring System of the Comprehensive Nuclear-Test Ban Treaty Organization (CTBTO), the responsibility of the On Site Inspection (OSI) division is the investigation of the source area to collect evidence that reveals whether a nuclear test has been conducted and, if the circumstances permit, to get a final localization of ground zero. At the location of an underground nuclear explosion, a damaged zone is expected to be present, including a cavity. Thus, cavity detection might become a major tool for the OSI division. 

In order to contribute to the method design, we investigate the scattering of the seismic wave field in the presence of an acoustic inclusion. The underlying technical questions of the OSI are still quite new and there are only few experimental examples that have been suitably documented to build a proper scientific groundwork. This motivates the investigation of the wave field on a purely numerical level and the simulation of potential observations based on recent advances in {numerical modeling of wave propagation problems.

As much as this is a challenging task in the applied fields, it is also interesting from a modeling and computational point of view. The classical scattering problem considers the wave propagation in an acoustic medium with an elastic obstacle, whereas we focus on the inverse situation of an elastic medium with an acoustic inclusion. 

For very simple cases the propagation of seismic waves can be described analytically \cite{Love1927treatise,Achenbach1973wave,aki1980quantitative,ben1981seismic}. For more complex cases, seismic waves with a significantly smaller or larger wave length than the characteristic size of the obstacle can be approximated by ray tracing methods \cite{cerveny2001seismic} or effective medium methods \cite{jordan2015effective,avila2016numerical,rahimzadeh2014comparative}, respectively. 
However, we are interested in the scattered wave patterns when the wavelengths of the propagating waves and the characteristic size of heterogeneities are comparable and here numerical methods become essential. There are many textbooks discussing the numerical modeling of seismic wave propagation  \cite{kelly1990numerical,carcione2001wave,cohen2001higher,moczo2014finite,igel2017computational,
bathe1996finite,zienkiewicz2000finite}. Here, SPEED \cite{mazzieri2013speed,Antonietti2014High,antonietti2012non,Faccioli1997elastic} is applied to a three-dimensional (3D) elastic-acoustic scattering problem for comparison. We consider in particular a 3D scattering problem consisting of a low-velocity
spherical acoustic inclusion embedded in a high-velocity elastic medium, whereby a
plane P-wave is scattered by the inclusion having a diameter similar to the P-wave’s wavelength. 
For this case Korneev and Johnson \cite{korneev1993scattering} provide an analytic solution which is used as the reference solution for the comparison of the numerical results. A similar study in 2D has been discussed in \cite{frehner2008comparison}.
Based on the analytic solution presented in \cite{korneev1993scattering}  
the investigation of seismic resonances origin to an acoustic inclusion was also discussed in \cite{schneider2017seismic}.

\section{Problem formulation}

Let $\Omega = \Omega_a \cup \Omega_e$ be an open bounded set $\Omega \subset \mathbb{R}^3$, having a spherical obstacle/acoustic inclusion $\Omega_a$ located at the origin as illustrated in Figure~\ref{fig:design}. We consider the elastic wave propagation problem in $\Omega$ described by
\begin{align}\label{eq:strong_problem}
\begin{cases}
\rho_e \ddot{\u}_e - \nabla \cdot \sig(\u_e) = {\bf f}, & \text{ in } \quad  \Omega_e,\\
\rho_a \ddot{\u}_a - \nabla \cdot \sig(\u_a) = {\bf 0}, & \text{ in } \quad  \Omega_a,\\
\text{ + coupling conditions}, & \text{ on }\quad \Gamma_I,\\
\text{ + boundary conditions}, & \text{ on } \quad\Gamma_B ,\\
\dot{\u}_e = \u_e = {\bf 0}, & \text{ in } \quad  \Omega_e,\\
\dot{\u}_a = \u_a = {\bf 0}, & \text{ in } \quad  \Omega_a,
\end{cases}
\end{align}
where $\rho_i$ is the mass density within the subdomain $\Omega_i$, $i=\{e,a\}$, and  $\u_i$ is the corresponding  displacement unknown vector. The Cauchy stress tensor $\sig:\mathbb{R}^{3\times 3}\rightarrow \mathbb{R}^{3\times 3}$ is expressed with the domain-wise constant Lamé parameters $\lambda_i$ and $\mu_i$ by
\begin{equation*}
\label{eq:stress}
\sig(\u_i) = \lambda_i(\nabla\cdot\u_i)\mathbf{I}+\mu_i (\nabla\u_i + \nabla\u_i^\top), \quad i =\{e,a\}.
\end{equation*}
\begin{figure}[h!]
\centering
\includegraphics[width=0.65\textwidth, trim={0cm 0pt 0cm 0pt} , clip=true]{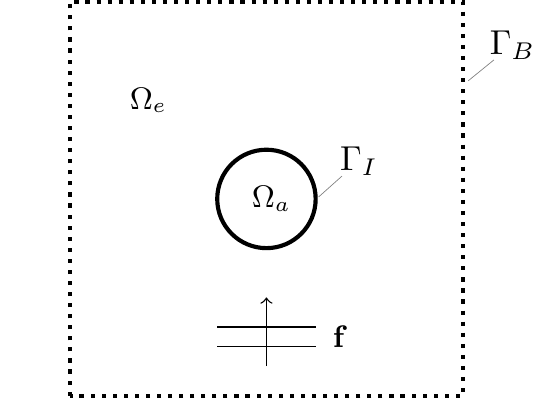}
\caption{Sketch of the domain $\Omega$. Homogeneous elastic medium $\Omega_e$ surrounding a spherical cavity $\Omega_a$. Force load ${\bf f}$ applied at the bottom of the domain.}
\label{fig:design}
\end{figure} 
We recall that the Lam\'e parameters are related to the pressure ($v_p$) and shear ($v_s$) velocity of the media as follows
\begin{equation}
\label{eq:lame}
\left\{\begin{array}{ll} \lambda_e = \rho_e(v_{p,e}^2-2v_{s,e}^2), &\quad \text{ in } \Omega_e, \\ \lambda_a = \rho_a v_{p,a}^2,  &\quad \text{ in } \Omega_a,\end{array}\right. 
\!\!\!\quad\quad
\left\{\begin{array}{ll} \mu_e = \rho_e v_{s,e}^2, & \quad\text{ in } \Omega_e, \\ \mu_a = 0, & \quad\text{ in } \Omega_a. \end{array}\right.
\end{equation}
Moreover, we also define the subdomain-wise constant wave number by
\begin{equation}
\label{eq:wavenumber}
\left\{\begin{array}{ll} k_{p,e} = \omega/v_{p,e}, & \quad \text{ in } \Omega_e, \\ k_{p,a} = \omega/v_{p,a}, &\quad \text{ in } \Omega_a,\end{array}\right. 
\!\!\!\quad\quad
\left\{\begin{array}{ll} k_{s,e} = \omega/v_{s,e}, & \quad \text{ in } \Omega_e, \\ k_{s,a} = 0, & \quad \text{ in } \Omega_a.\end{array}\right. 
\end{equation}
Since we are interested in the scattered wave field subject to an incoming pressure plane wave, we consider a distribution of body force given by
\begin{equation}
\f(\x,t) = \phi(t)\delta(z-z_0)\bm{e}_3,
\label{eq:rhs}
\end{equation}
with time profile $\phi(t)$ acting on the plane $z = z_0$ for some $z_0$ far from the acoustic inclusion at the bottom of the domain. 

\subsection{Coupling conditions}
\label{sec:inter}
At the interface between $\Omega_e$ and $\Omega_a$ we request only continuity in the normal component of the displacement as well as of the stress tensor:
\begin{align}
\u_e \cdot\n_a &= \u_a\cdot\n_a, \label{eq:u_jump}\\
\sigma(\u_e)\n_a &= \sigma(\u_a)\n_a,\label{eq:s_jump}
\end{align}
where $\n_a$ is the normal vector pointing outward from $\Omega_a$.
Note that equation~\eqref{eq:u_jump} is equivalent to impose a null jump for the normal component of the displacement field across the interface between the elastic and acoustic domain, that is $[[\u\cdot\n_a]]=\bm{0}$ on $\Gamma_I$. However, if the difference of displacement in tangential direction stays small, this condition can be replaced by the jump condition $[[\u]] = \bm{0}$, which will be used in the DG formulation proposed in Section~\ref{sec:dgform}. Note also that equation~\eqref{eq:s_jump} includes a free surface condition in the tangential components for elastic medium.

\subsection{Boundary conditions}\label{sec:absorbing}
A possible approach to approximate the radiation condition for the unbounded domain $\mathbb{R}^3$ consists in modeling an absorbing boundary layer by the introduction of a fictitious traction ${\bf t^*}$
on $\Gamma_B$.  
Here, we consider the local P3 paraxial conditions presented in \cite{stacey1988improved}, which is sufficiently accurate if $v_{p,e} /v_{s,e} \leq 2$, as in the application under consideration.
More specifically, the P3 paraxial absorbing conditions read as
\begin{align*}
\partial_{n_e}(\u_e\cdot\n_e) &= -\frac{1}{v_{p,e}}\partial_t(\u_e\cdot\n) - \frac{v_{p,e}-v_{s,e}}{v_{p,e}}\left[ \partial_{\tau_1}(\u_e\cdot\bm{\tau_1}) + \partial_{\tau_2}(\u_e\cdot\bm{\tau_2})\right]\\
\partial_{n_e}(\u_e\cdot\bm{\tau_1}) & = -\frac{1}{v_{s,e}}\partial_t(\u_e\cdot\bm{\tau_1}) - \frac{v_{p,e}-v_{s,e}}{v_{p,e}} \partial_{\tau_1}(\u_e\cdot\n) \\
\partial_{n_e}(\u_e\cdot\bm{\tau_2}) & = -\frac{1}{v_{s,e}}\partial_t(\u_e\cdot\bm{\tau_2}) - \frac{v_{p,e}-v_{s,e}}{v_{p,e}} \partial_{\tau_2}(\u_e\cdot\n) 
\end{align*}
where $\bm{\tau_1}$ and $\bm{\tau_2}$ are two mutually orthogonal unit normal vectors on the
plane orthogonal to normal vector $\n_e$ pointing outward of $\Omega_e$.  $\bm{\tau_1}$ and $\bm{\tau_2}$ span the tangent plane to the surface $\Gamma_B$ in each point such that $\{\bm{\tau_1},\bm{\tau_2},\n_e\}$ is a right handed Cartesian system. The traction term ${\bf t}^* = \sig^*(\u_e)\n_e$ defined on the absorbing boundary in the local coordinate system $(\bm{\tau_1}, \bm{\tau_2}, \n_e)$ has then the following expression
\begin{equation*}
\left[ \begin{matrix}
t^*_{\bm{\tau_1}} \\
t^*_{\bm{\tau_2}} \\
t^*_{\n_e}
\end{matrix}
\right] =  \left[ \begin{matrix}
\frac{\mu_e(2v_{p,e}-v_{s,e})}{v_{s,e}}\partial_{\tau_1}(\u_e\cdot \n_e) - \frac{\mu_e}{v_{s,e}}\partial_{t}(\u_e\cdot\bm{\tau_1})\\
\frac{\mu_e(2v_{p,e}-v_{s,e})}{v_{s,e}}\partial_{\tau_2}(\u_e\cdot \n_e) - \frac{\mu_e}{v_{s,e}}\partial_{t}(\u_e\cdot\bm{\tau_2})\\
\frac{\lambda_e v_{s,e} + 2\mu_e(v_{p,e}-v_{s,e})}{v_{s,e}}\left[\partial_{\bm{\tau_2}}(\u_e\cdot\bm{\tau_1}) + \partial_{\bm{\tau_1}}(\u_e\cdot\bm{\tau_2})\right] - \frac{\lambda_e+2\mu_e}{v_{s,e}}\partial_{t}(\u_e\cdot\n_e)
\end{matrix}
\right],
\end{equation*}
that can be easily rewritten in term of the gobal coordinate system $(x,y,z)$. See \cite{casadei1998implementation} for more details.

\section{Numerical discretization}
\label{sec:num}

Piece-wise constant material parameters result in contrasting wave lengths and give reason to approximate the solution with distinct discretization parameters in each domain $\Omega_i$, $i = \{e,a\}$. Especially when the velocity contrast is comparably high, this motivates to use proper space discretization parameters in each subdomain, in order to catch the main features of the wave phenomenon. This motivates the choice of the following Discontinuous Galerkin numerical discretization.

\subsection{Mesh and trace operators}

We consider a (not necessarily conforming) decomposition $\T_\Omega$ of $\Omega$ into two nonoverlapping polyhedral sub-domains $\Omega_e$ and $\Omega_a$, i.e., $\bar{\Omega} =  \Omega_e \cup \bar{\Omega}_a$, $\Omega_e \cap \Omega_{a} = \emptyset$. On each $\Omega_i$, $i=\{e,a\}$ we consider a conforming, quasi-uniform computational mesh $\T_{h_i}$ of granularity $h_i > 0$ made by open disjoint elements $K_i^j$, and suppose that each $K_i^j \in \Omega_i$ is the image through a bilinear map $\Phi_i^j: \hat{K} \to K_i^j$ of the  reference hexahedron $\hat{K}=[-1,1]^2$.
We define an interior face $F$ as the non-empty interior of $\partial K_e \cap \partial K_a$, for some $K_e \in \T_{h_{e}}$ and $K_a \in \T_{h_{a}}$, and collect all the interior faces in the set $F_h^I$. Moreover, we
define $F_h^{B}$ as the sets of all boundary faces where absorbing boundary conditions are imposed. 
Finally, we assume that for any element $K \in \T_h$ and for any face $F\subset \partial K$ it holds $h_K \lesssim h_F$. For more details see \cite{GeorgoulisHallHouston2007b,PerugiaSchotzau_2002} as well as  \cite{Dryja_2003,DryjaSarkis_2007} for the case of highly discontinuous coefficients.

Let $K_e \in \T_{h_{e}}$ and $K_a \in \T_{h_{a}}$ be two elements sharing a face $F \in F_h^I$, and let $\n_i$ be the unit normal vectors to $F$ pointing outward to $K_i, i\in\{e,a\}$, respectively. For (regular enough) vector and tensor-valued functions $\v$ and $\boldsymbol{\tau}$,
we denote by $\v_i$ and $\boldsymbol{\tau}_i$ the traces of $\v$ and $\boldsymbol{\tau}$ on $F$, taken within the interior of $K_i, i\in\{e,a\}$, respectively, and set
\begin{equation*}
[[ \v ]] = \v_e \odot \n_e + \v_a \odot \n_a, \quad 
[[ \boldsymbol{\tau} ]] = \boldsymbol{\tau}_e \n_e + \boldsymbol{\tau}_a \n_a, \quad 
\{ \v \} = \dfrac{\v_e + \v_a}{2}, \quad 
\{ \boldsymbol{\tau} \} = \dfrac{\boldsymbol{\tau}_e + \boldsymbol{\tau}_a}{2},
\end{equation*}
where $\v \odot \n = (\v^T\n+\n^T\v)/2$. On $F \in F_h^B$, we set $\{\v\} = \v$, $\{\boldsymbol{\tau} \} = \boldsymbol{\tau}$, $[[ \v]] = \v \odot \n$, $[[\boldsymbol{\tau} ]] = \boldsymbol{\tau} \n$.

\subsection{Discontinuous Galerkin Sprectral Element discretization}\label{sec:dgform}

For each subdomain $\Omega_i$, $i=\{e,a\}$ we consider a nonnegative integer $N_i$, and we define the finite dimensional space
\begin{equation*}
V_{h_i}^{N_i}(\Omega_i) = \{ \v \in {\bf C}^0(\bar{\Omega}_i) : (\v_{|_{K_i^j}} \circ \Phi^j_i) \in [ \mathbb{P}^{N_i}(\hat{K})]^3 \quad \forall K^j_i \in \T_{h_i} \},
\end{equation*}
where $\mathbb{P}^{N_i}(\hat{K})$ is the space of polynomials of degree $N_i$ in each coordinate direction on $\hat{K}$. Then, we define the finite dimensional trial space $\VDG$ as $\VDG = \prod_{i=e,a} V_{h_i}^{N_i}(\Omega_i)$. The semidiscrete Discontinuous Galerkin approximation of problem \eqref{eq:strong_problem} reads: $\forall t \in (0,T]$, find $\u_h = \u_h(t) \in \VDG$ such that
\begin{equation}\label{eq:dgformulation}
\sum_{i=e,a}\int_{\Omega_i} \rho_i\ddot{\u}_h(t) \cdot \v \tin\Omega + \mathcal{A}_h(\u_h(t),\v) = \mathcal{F}_h(\v)
\qquad \forall \v \in \VDG,
\end{equation}
subjected to the initial conditions $\dot{\u}_h(0) = \u_h(0) = {\bf 0}$. 
The right hand side $\mathcal{F}_h(\cdot)$ is defined as
\begin{equation*}
\mathcal{F}_h(\v)  = \int_{\Omega_e} {\bf f}(t) \cdot \v \tin\Omega +\int_{\Gamma_B} {\bf t}^* \cdot \v \tin \Gamma \qquad \forall \v \in \VDG,
\end{equation*}
while the bilinear form $\mathcal{A}_h(\cdot,\cdot) $ as
\begin{multline}\label{def:bilinearform}
\mathcal{A}_h(\u,\v) = \sum_{i=e,a}\int_{\Omega_i} \sig(\u) : \bm{\epsilon}(\v) \tin\Omega \\ 
- \sum_{F\in F_h^I}\Big( \int_F \{ \sig(\u) \}:[[\v]]\tin\Gamma  +  \int_F \{ \sig(\v) \}:[[\u]] \tin\Gamma  - \int_F \eta[[\u]]:[[\v]] \tin\Gamma \Big) ,    
\end{multline}
for any $\u,\v \in \VDG$, being $\bm{\epsilon}(\v) = (\nabla \v + \nabla \v^T)/2$ and $\eta$ a positive parameter to be choosen large enough, cf. \cite{antonietti2012non}.

\begin{remark}
Implicit in the derivation of formulation \eqref{eq:dgformulation} is the use of coupling conditions described in section \ref{sec:inter}. For the sake of presentation we derive formulation \eqref{eq:dgformulation} in the case of a partition made by two subdomain $\Omega_e$ and $\Omega_a$. However, it can be easily extended for accommodating  different elastic or acoustic subdomains, as it will be considered in Section \ref{sec:results}. 
Finally, note that the discrete solution is piecewise discontinuous across macro elements $\Omega_i, i\in \{e,a\}$ and (weak) continuity is enforced based on employing, at a subdomain level, the symmetric interior penalty DG (SIPG) method \cite{antonietti2012non}. We refer to \cite{antonietti2016stability} for a unified analysis of the $h$-version of the method and to \cite{AntoniettiFerroniMazzieriQuarteroni_2016b,antonietti2015high} for the $hp-$version of the method and it analysis.
\end{remark}

\subsection{Fully discrete formulation}

In this section we present the time integration of the semi-discrete formulation \eqref{eq:dgformulation}. By fixing a basis for the discrete space $\VDG$, the semi-discrete algebraic formulation of problem \eqref{eq:dgformulation}, reads as
\begin{equation}
\bm{M}_0\ddot{\U}(t) +  \bm{M}_1\dot{\U}(t) + (\bm{M}_2 + \bm{A}){\U}(t) = \bm{F}(t) \qquad \forall t \in (0,T],
\label{eq:semi}
\end{equation}
supplemented by the initial conditions $\dot{\U}(0) = \U(0) = \bm{0}$. Here, denoting by $N_{dof}$ the total number of degrees of freedom, the vector $\U = \U(t) \in \mathbb{R}^{N_{dof}}$ contains, for any time $t$, the expansion coefficients of the semi-discret solution $\u_h(t) \in \VDG$ in the chosen set of basis functions. Analogoulsy, $\bm{M}_0$ and $\bm{A}$ are the matrices representations of the bilinear forms 
\begin{equation*}
\sum_{i=e,a}\int_{\Omega_i} \rho_i\ddot{\u}_h(t) \cdot \v \tin\Omega \; \; \text{ and } \;\; \mathcal{A}_h(\u_h(t),\v),
\end{equation*}
respectively, cf. \eqref{eq:dgformulation}. Inserting the absorbing conditions from Section \ref{sec:absorbing} for the boundary term $\int_{\Gamma_B} \bm{t}^* \cdot \v \tin \Gamma$ give rise to the matrices $\bm{M}_1$ and $\bm{M}_2$ in equation \eqref{eq:semi}. Finally $\bm{F}$ is the vector representation of the linear functional $\mathcal{F}_h(\cdot)$ containing the body force term $\f$, cf. \eqref{eq:dgformulation}.

For the time integration of the system of second order ordinary differential equations \eqref{eq:semi}, we employ the leap-frog method \cite{quarteroni2008numerical}, that is a  widely employed time marching scheme for the numerical simulation of elastic waves propagation, see for example \cite{bao1998large,chaljub2007spectral,komatitsch1999introduction,mercerat2015nodal}.
With this aim we subdivide the time interval $(0,T]$ into $N_T$ subintervals of amplitude $\Delta t$, and we denote by $\U_{n} \approx  \U(t_n)$ and $\bm{F}_{n} \approx  \bm{F}(t_n)$  the approximation of $\U$ and $\bm{F}$ at time $t_i=i\Delta t$, $i=1,2,...,N_T$, respectively.
System \eqref{eq:semi} approximated with the leap-frog scheme reads as:
\begin{align}
\bm{M}_0 \U_1 & =  \frac{\Delta t^2}{2}\bm{F}_0, \label{leap-frog1} \\
 (\bm{M}_0 + \frac{\Delta t}{2}\bm{M}_1) \U_{n+1} & =  (2\bm{M}_0-\Delta t^2 \bm{Q})\U_n + (\bm{M}_0 - \frac{\Delta t}{2}\bm{M}_2) \U_{n-1} + \Delta t^2 \bm{F}_n,   \label{leap-frog2}
\end{align}
for  $n= 1,...,N_T-1$, with $\bm{Q} = \bm{A} + \bm{M}_2$. 
We notice that \eqref{leap-frog2} involves a linear system with  matrix  $\bm{M}_0 - \frac{\Delta t}{2}\bm{M}_2$ to be solved at each time step. 
The choice of the  basis functions spanning the space $\VDG$ strongly influences the structure of the  matrix $\bm{M}_0 - \frac{\Delta t}{2}\bm{M}_2$ and, therefore,
the computational cost related to the solution of the linear system. Furthermore, since the leap-frog method is an explicit second order accurate scheme, to ensure its numerical stability a Courant - Friedrich - Levy (CFL) condition has to be satisfied (see \cite{quarteroni2008numerical}).

\section{Analytic solution}

For the following analysis, we consider the total wave field $\u$ expressed as the sum of the incident ($\u_I$) and the scattered ($\u_S$) wave fields as follow
\begin{equation}
\u(\x,t) = \u_I(\x,t)+\u_S(\x,t).
\label{eq:split}
\end{equation}
In a homogeneous elastic domain the incident wave origin to a body load \eqref{eq:rhs} has only a contribution in the 3rd component which is given by
\begin{equation}
u_{I}^3(\x,t) = \frac{1}{2\rho_e v_{p,e}} H(t-\frac{|z-z_0|}{v_{p,e}})
\int_0^{t-\frac{|z-z_0|}{v_{p,e}}} \phi(\tau)\tin\tau.
\label{eq:travel_plane_1}
\end{equation}
Note that $\u_I$ solves then the elastic wave question $\rho_e \ddot{\u}_I + \nabla\cdot \sig(\u_I) = \f$ in $\mathbb{R}^3$ with parameters from $\Omega_e$.  
Hence, to recover the total wave field it is sufficient to solve the following problem for $\u_S$:  
\begin{align}\label{eq:scat_th}
\begin{cases}
\rho_e \ddot{\u}_S - \nabla \cdot \sig(\u_S) = {\bf 0}, & \text{ in } \quad  \Omega_e,\\
\rho_a \ddot{\u}_S - \nabla \cdot \sig(\u_S) = -\rho_a \ddot{\u}_I + \nabla\cdot \sig(\u_I), & \text{ in } \quad  \Omega_a,
\end{cases}
\end{align}
together with solely absorbing conditions at the boundary $\Gamma_B$, cf. Figure~\ref{fig:design}. This approach has also been used in \cite{esterhazy2017application}.

In case of a time-harmonic force load outside the domain (at infinity), the body force equals zero, i.e. $\f=\bm{0}$. For this case, an analytic solution is presented by Korneev and Johnson\cite{korneev1993scattering} and summarized here for the sake of completeness. 
In the time-harmonic case, the incident as well as the scattered wave fields can be expressed as
\begin{equation}
\u_I(\x,t) = \mathfrak{R}\{\U_I(\x)e^{-i\omega t}\}, \quad \u_S(\x,t) = \mathfrak{R}\{\U_S(\x)e^{-i\omega t}\}
\label{timeharmonic}
\end{equation}
where $\U_I(\x)$ and $\U_S(\x)$ are complex-valued functions. In particular, the interaction of the incident wave with the sphere gives rise to a scattered displacement field inside as well as outside of the sphere. To this end we omit the time dependence and use the following notation
\begin{equation}
\U_{1} (\x) = \U_S(\x)|_{\Omega_a},\quad \U_{2} (\x) = \U_S(\x)|_{\Omega_e}.
\end{equation}
In order to construct an analytical solution \cite{korneev1993scattering} used the system of spherical vectors in the spherical coordinte system $(r,\theta,\phi)$ with unit vectors $\{  \bm{\hat{r}}, \bm{\hat{\theta}},\bm{\hat{\phi}}\}$ developed by Petrashen \cite{Petrashen1945}
\begin{align*}
\Y_{lm}^0 = \Y_{lm}^0 (r,\theta,\phi) &= \mathbf{r}\times\nabla Y_{lm}\\
\Y_{lm}^+ = \Y_{lm}^+ (r,\theta,\phi)   &= (l+1)\mathbf{\hat r}Y_{lm}-r\nabla Y_{lm}\\
\Y_{lm}^- = \Y_{lm}^- (r,\theta,\phi)   &= l\mathbf{\hat r}Y_{lm}+r\nabla Y_{lm}
\end{align*}
where $r$ is distance from the center of the sphere with $\bm{r} = r\bm{\hat{r}}$
and $Y_{lm}$ are the unnormalized spherical harmonic functions, defined as
\begin{equation*}
Y_{lm}=Y_{lm}(\theta,\phi) = e^{im\phi}P_l^m(\cos(\theta))
\end{equation*}
such that an arbitrary vector function $\U$ can be represented in the form
\begin{equation*}
\U(\x) = \sum_{l,m} a_{lm}^0(r)\Y_{lm}^0+a_{lm}^+(r)\Y_{lm}^++a_{lm}^-(r)\Y_{lm}^-.
\end{equation*}
In this coordinate system  an incident plane harmonic P-wave, propagating in the positive z-direction in $\Omega_e$ is given by
\begin{equation*}
\U_I^P(\x) 
= \sum \left\{j_{l+1}(k_{2,p}r)\Y_{l0}^+ + j_{l-1}(k_{2,p}r)\Y_{l0}^-\right\} e^{-\frac{i\pi}{2}(l+1)}
\end{equation*}
where the $j_l(z)$ are the spherical Bessel functions. Furthermore it is possible to express the scattered wave field inside and outside the sphere separately by
\begin{align*}
 \U_1=&\sum_{l\geq 0}\Big\{\Big({a_l^{(1)}}j_{l+1}(k_{p,1}r)+ l{b_l^{(1)}}j_{l+1}(k_{s,1}r) \Big) \Y^{+}_{ l0}  \\
&+\Big(-{a_l^{(1)}}j_{l-1}(k_{p,1}r) + (l+1) {b_l^{(1)}}j_{l-1}(k_{s,1}r)\Big) \Y^{-}_{l0} \Big\}e^{- \frac{i \pi }{2}(l+1)},
\end{align*}
\begin{align*}
 \U_2=&\sum_{l\geq 0}\Big\{\Big({a_l^{(2)}}h_{l+1}(k_{p,2}r)+ l{b_l^{(2)}}h_{l+1}(k_{s,2}r) \Big) \Y^{+}_{ l0}  \\
&+\Big(-{a_l^{(2)}}j_{l-1}(k_{p,2}r) + (l+1) {b_l^{(2)}}j_{l-1}(k_{s,2}r)\Big) \Y^{-}_{l0} \Big\}e^{- \frac{i \pi }{2}(l+1)}, 
\end{align*}
respectively, where $h_j(z)$ are the spherical Hankel functions of second kind and
$a_l^{(\nu)},b_l^{(\nu)}, \nu=1,2$ are the coefficients which are given
explicitly in \cite{korneev1993scattering}. The wavenumbers $k_{p,i}$ and $k_{s,i}$ are given by $\omega/v_{p,i}$ and $\omega/v_{s,i}$, where $v_{p,i}$ and $v_{s,i}$ are the propagation velocities of P- and S-waves inside ($i=a$) and outside ($i=e$) of the cavity, respectively.
The unknown coefficients $a_l^{(\nu)},b_l^{(\nu)}$ can be determined by solving a linear system that arises from the following continuity conditions, which are valid at the acoustic-elastic interface:
\begin{equation*}
\U_1 \cdot \n = (\U^P_I+\U_2)\cdot \n \quad \text{ and } \quad \bm{\sigma}(\U_1) \, \n = \bm{\sigma}(\U_I^P+\U_2) \, \n
\end{equation*}
The first condition describes the continuity of the normal component of the
displacement. The tangential displacement components  are free due to the
fact that no shear stress can be  transmitted to the acoustic domain.  Since $\mu=0$ in $\Omega_a$
the traction vector $\bm{\sigma}(\U_1) \n$ points in the normal direction $\n$ with respect to the interface.
Thus the second condition forces the tangential components of the traction
vector $\bm{\sigma}(\U_I^P+\U_2)\n\cdot \mathbf{t}$ to be zero, for any vector
$\mathbf{t}$ which is orthogonal to $\n$. Thus, the acoustic-elastic interface
acts as a free surface for the components tangential to the interface and
transmits only normal components of displacement and stress between  the acoustic
and elastic domains.

The total wave-field outside and inside the cavity is given by 
\begin{equation*}
\U_{tot} = \U^P_I+\U_2 \quad \text{ and } \quad \U_{tot} = \U_1, \quad \text{ respectively.}
\end{equation*}
In order to retrieve the solution of the wave equation for a time dependent incident field,
the scattering problem is solved for many frequencies and 
the time harmonic functions (Eq. \ref{timeharmonic}) are combined by applying inverse Fourier transform.
In order to compare the analytical with numerical solutions  all wave fields are convolved with the Ricker wavelet
\begin{equation*}
R(t)=(1-2\pi^2 f_0^2 t^2) e^{\pi^2 f_0^2 t^2}, 
\end{equation*}
describing a common seismic model wavelet. 
However, any other wavelet can be used instead in order to retrieve arbitrary time histories.
With the convolution theorem the calculated synthetic seismogram can be expressed as
\begin{equation*}
s(t) = \mathfrak{F}^{-1}\Bigl\lbrack \; \U(\{\omega\}) \; |\mathfrak{F}[R(t)](\{\omega\}) |\; \Bigr\rbrack (t),
\end{equation*}
where $\mathfrak{F}$ and $\mathfrak{F}^{-1}$ are  Fourier's transform and it's discrete inverse, respectively. $\{\omega\}$ is the  set of frequencies for which the solution is computated.

\section{Results} 
\label{sec:results}

In this section we want to address a 3D scattering wave propagation problem consisting of a low-velocity
spherical acoustic inclusion embedded in a high-velocity elastic medium.  In particular we want to compare our numerical results with respect to the analytical one provided by  Korneev and Johnson in \cite{korneev1993scattering}.

\subsection{Mesh generation}

Special attention must be given to the grid generation as meshing a spherical inclusion inside a cube with hexahedrons is not a trivial task. Especially as in this case, when the wave length inside the sphere is much smaller than outside. This gives reason to chose a smaller mesh size inside the inclusion. 
Using non-curved elements non-conforming meshes inside and outside the sphere will lead to empty and overlapping regions which would lead to numerical instabilities and must therefore be avoided. As a work-around we added another small box around the sphere such that the non-conforming interface can be generated between the small and the big cube while having a conforming interface between the small cube and the sphere, see Figure \ref{fig:mesh}. A strategy to overcome this issue is presented in \cite{RODRIGUEZROZAS201644}.
However, the parameters are discontinuous across the boundary of the sphere and therefore we can select different discretization parameters inside and outside the spherical cavity. 
In summary, DG jumps are applied to both interfaces: the non-conforming interface between the two elastic cubes and the conforming interface between the acoustic and elastic domain where physical parameters are discontinuous.

\begin{figure}[h]
\centering
\includegraphics[width=0.75\textwidth]{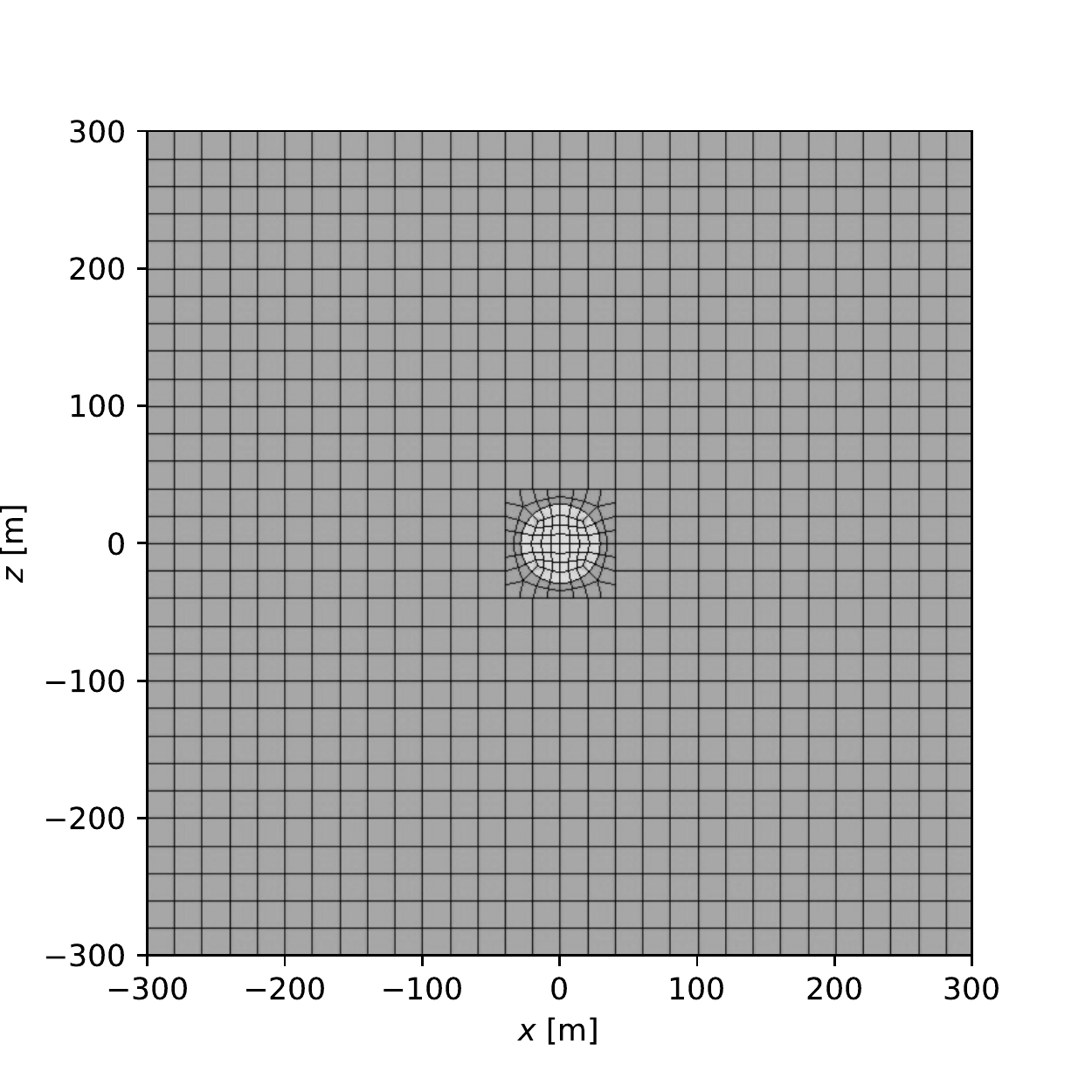}
\caption{Illustration of the mesh strategy using Trelis (\texttt{http://www.csimsoft.com/trelis.jsp}). The dark and light gray domain correspond to the elastic and acoustic cavity, respectively. The spherical domain is connected with conforming interface to a surrounding smaller box with the same small mesh size. The smaller box is embedded in a bigger box with a larger mesh size resulting in a non-conforming interface.}
\label{fig:mesh}
\end{figure}

\subsection{Input data}
\label{sec:input}
 The mechanical parameters for the elastic and acoustic materials are given in Table~\ref{tab:input}.
\begin{table}[h]
\centering
\renewcommand{\arraystretch}{1.2}
\begin{tabular}{|c|c|c|c|}
\hline
Domain &  $\rho$ [kg/m$^3$]& $v_s$ [m/s] & $v_p$ [m/s]\\
\hline
$\Omega_a$ & 1000 & 0 & 1500\\
\hline
 $\Omega_e$ & 2700 & 2310 & 4000\\
 \hline
\end{tabular}
\caption{Physical parameters for the test case considered. }
\label{tab:input}
\end{table}

The time profile of our seismic source is described by the Ricker wavelet 
\begin{equation*}
R(t) = \left(1-2\beta (t-t_0)^2\right)e^{-\beta (t-t_0)^2}, \quad \beta = \left(\frac{\omega_p}{2}\right)^2 
\end{equation*}
where $\omega_p = 2\pi f_p$ is the angular peak frequency of the Ricker wave and $t_0$ a time offset. 
The shape of the Ricker profile is shown in Figure~\ref{fig:ricker} and more details on the frequency band can be found in \cite{Wang2015ricker}.

\begin{figure}[h]
\centering
\includegraphics[width=0.8\textwidth]{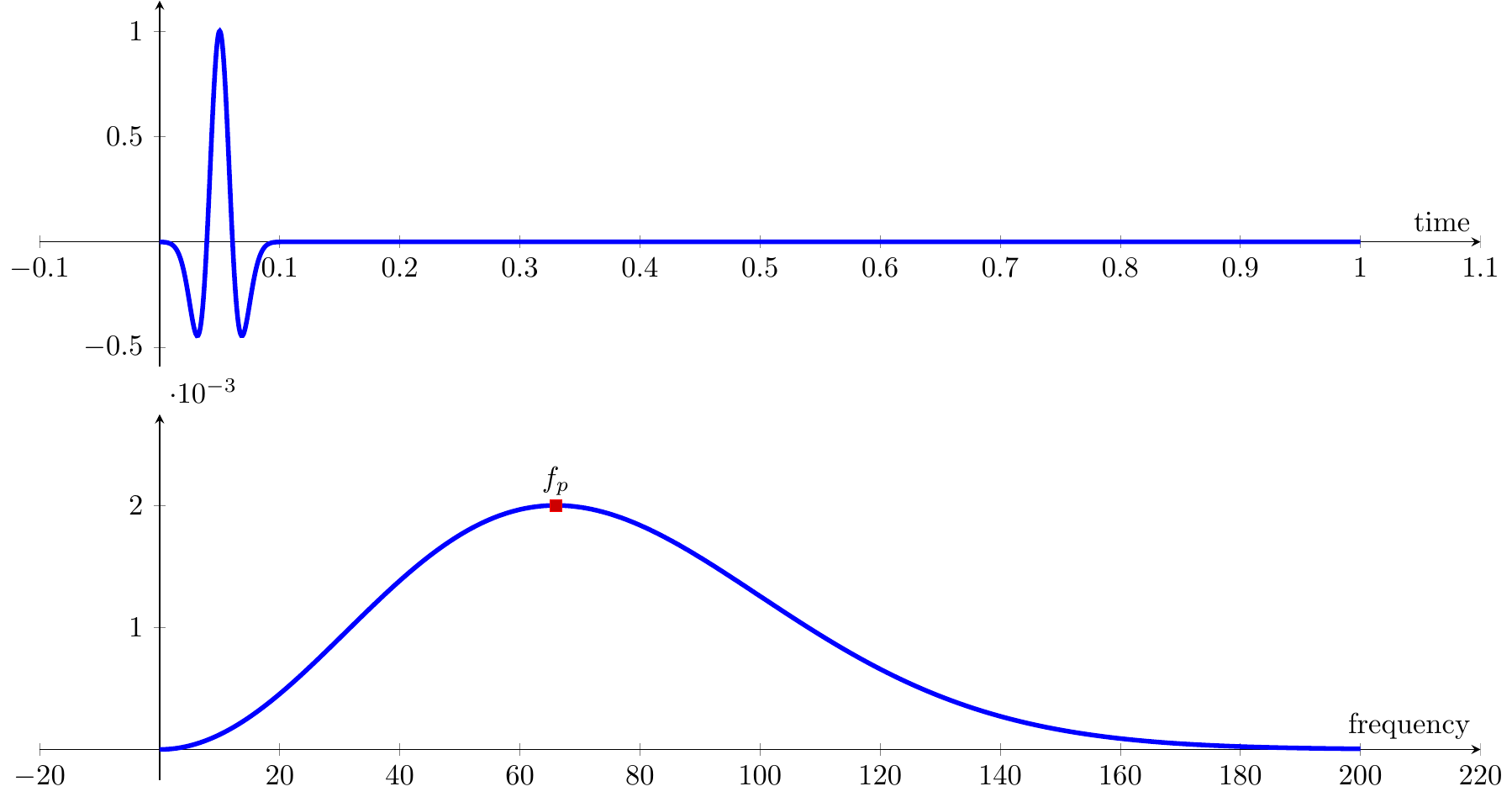}
\caption{Top: Time profile for the Ricker wave with peak frequency ${f_p= 66 ~\mathrm{[Hz]}}$ and time shift $t_0=0.03~[s]$. Bottom: Spectrum of the Ricker wave}
\label{fig:ricker}
\end{figure}

Furthermore we were interested in  the case of an incident Ricker wave pulse of a mean wave length of the size of the cavity, i.e. $\lambda = 2R$. Hence we chose:
\begin{align*}
\beta = 44000 \quad \Rightarrow \quad &f_{peak} \sim 66.7 Hz\quad \Rightarrow \lambda_{p,e} = v_{p,e}/f_{peak} \sim 60 = 2R,
\end{align*}
However, to minimize numerical dispersion errors we had to choose the grid size according to S-velocity and the maximum frequency in the Ricker spectrum as follows
\begin{align*}
f_{max} = 3f_{peak} \sim 200 Hz\quad \Rightarrow \quad \lambda_{s,e} = v_{s,e}/f_{max} \sim 11.5 \quad \Rightarrow \quad N_e=4, h_e = 5,
\end{align*}
in order to obtain the resolution of 10 points per wavelength (reasonable for small wave numbers) and avoid spurious numerical effect from the artificial boundary. 
With the same argumentation one should choose the grid size in the acoustic domain by 
\begin{equation*}
\lambda_{p,a} = v_{p,a}/f_{max} \sim 7.5 \quad \Rightarrow \quad N_a = 4, h_a = 3.5.
\end{equation*}
Note that we use different grid sizes in the domain resulting in a non-conforming mesh, but we use the same polynomial degree.
In our domain of interest $(600\times600\times600)~m^3$ this would lead to more than $17.e+6$ grid points and more than $1.e+9$ spectral nodes. The minimal grid size in this mesh is about $h_{min} = 1.1~m$ resulting in a time step size
\begin{equation*}
\Delta t = 0.2\times \frac{0.175\ h_{min}}{v_{p,max}} \sim 8.e{-6}~s.
\end{equation*}
For a simulation time until $T = 1~s$ this gives 125000 time steps.

However, due to limited computational resources we rather resolve for the peak frequency and put the boundary further away, i.e.
\begin{equation*}
\lambda_{s,e} = v_{s,e}/f_{peak} \sim 35 \text{ and } \lambda_{p,a} = v_{p,a}/f_{peak} \sim 22  \quad \Rightarrow \quad N_e = 4, h_e = 17, N_a = 10,
\end{equation*}
on a domain of dimension $(4000\times4000\times 2400)~m^3$ reducing the number of grid points to about 56000. With a minimal grid size $h_{min} = 4.74$ the corresponding time step size is $\Delta t = 4.e-5~s$ and 25000 time steps.

\subsection{Analysis of the results}

In order to provide an overview we show in Figure~\ref{fig:snap} snapshots of
the wave field in the XZ-plane. The incident plane P-wave travels with constant amplitude  through the
elastic domain from the bottom to the top with a velocity of $4000~m/s$ and reaches the
acoustic-elastic interface at about $0.1~s$.  With a positive impedance contrast
from the elastic to the acoustic domain given by the material parameters in
Table~\ref{tab:input} about 75\% of the incident wave is reflected resulting in
the primary scattered spherical P- and S- waves which can be seen at $t=0.15$~s.
About 25\% are transmitted into the cavity where it propagates only as
P-wave with a lower velocity of $1500~m/s$.  Each time the acoustic wave hits the
boundary of the cavity,  about 75\% of its energy are now reflected back and only
about 25\% are transmitted to the elastic domain. This yields an acoustic wave
energy trapped inside the cavity expressed in multiple reverberations that
couple out into the elastic medium periodically which can be seen for $t=0.2~s$ to $t=0.45~s$.
In \cite{schneider2017seismic} this resonance phenomenon has been investigated in more detail based on an analytic solution.

For validation we compare the analytic and the numerical solution along four profiles which are illustrated in Figure~\ref{fig:valid}. Profile A is a vertical profile for $x = 0$ and  the
$z-$coordinate ranging from $-100~m$ to $100~m$ with a distance of $2~m$ while profile
B is a horizontal profile for $z = 0$ and  the $x-$coordinate
ranging from $-100~m$ to $100~m$. The horizontal profiles C and D are located
further away  with vertical locations
at $300~m$ and $-300~m$ from the cavity, which is located at the origin.
Profile C and D range horizontally from $x=-300~m$ to $x=300~m$ with a distance of
$10~m$.

\begin{figure}[h!]
\centering
\includegraphics[width=0.8\textwidth, trim= 0cm 0cm 0cm 0cm, clip=true]{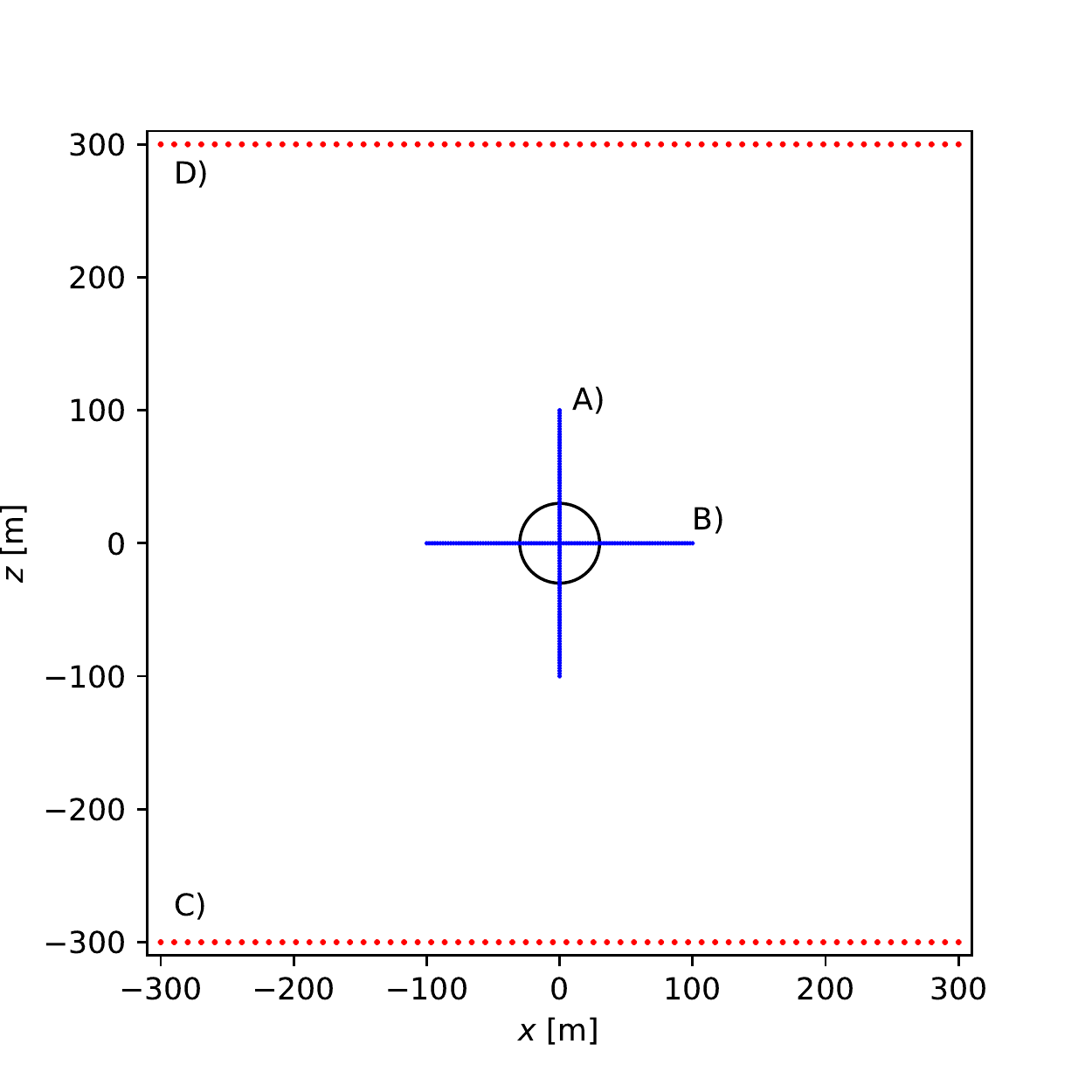}
\caption{Cross-section of the computational domain. Validation points along four profiles A),B),C) and D) are also represented.}
\label{fig:valid}
\end{figure}

\begin{figure}
\includegraphics[width=0.4\textwidth,scale=0.1]{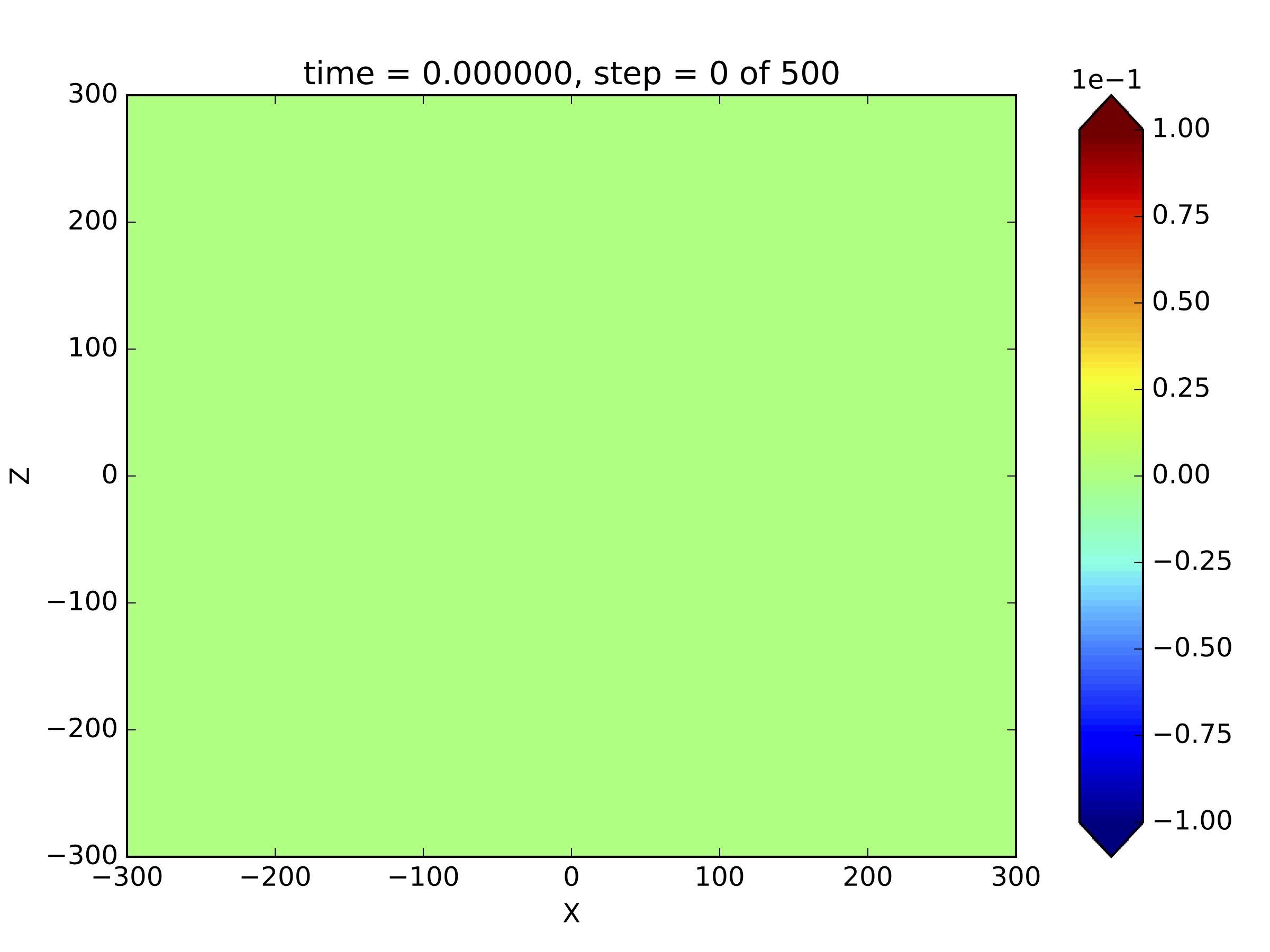}
\includegraphics[width=0.4\textwidth,scale=0.1]{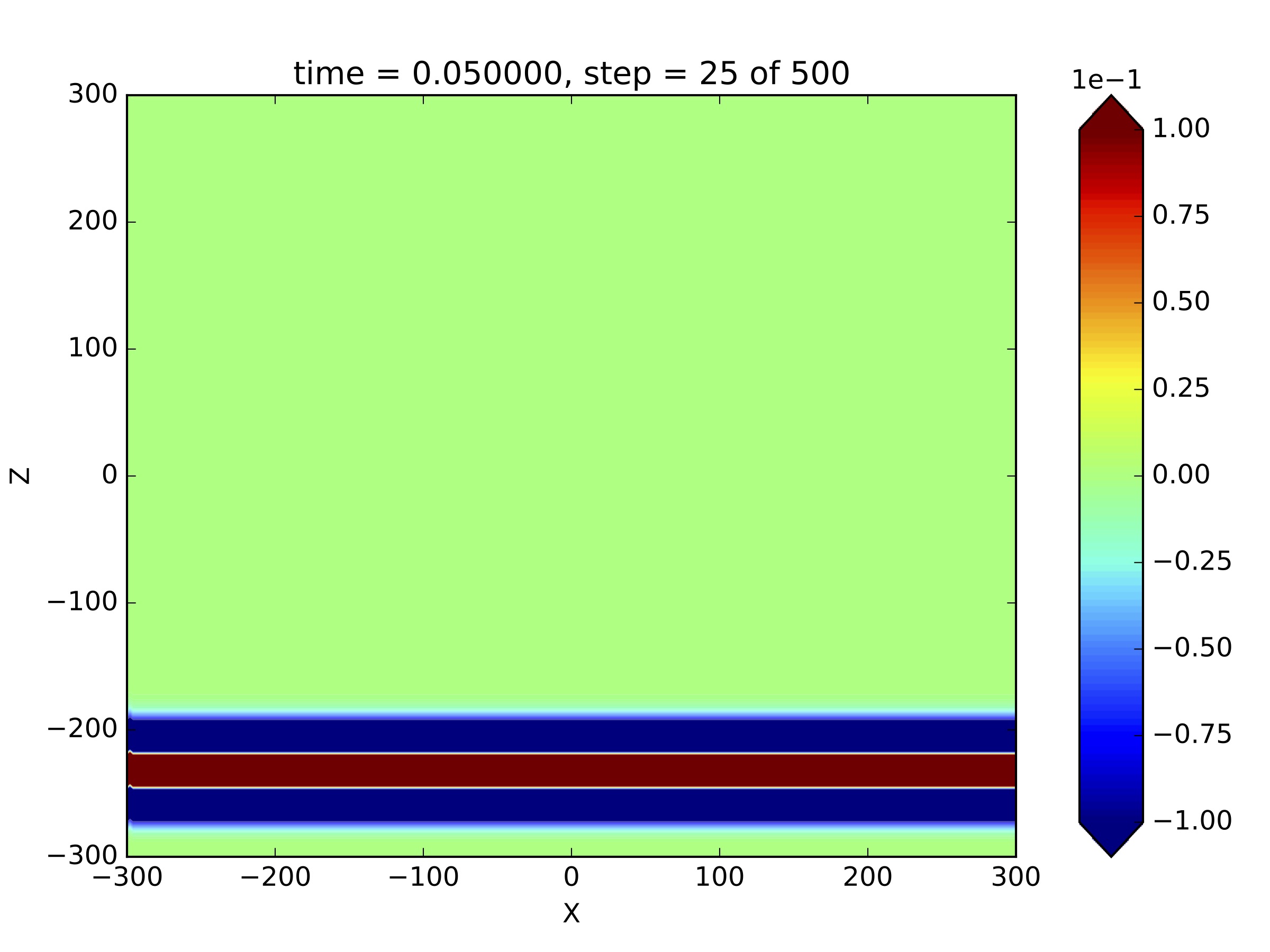}
\includegraphics[width=0.4\textwidth,scale=0.1]{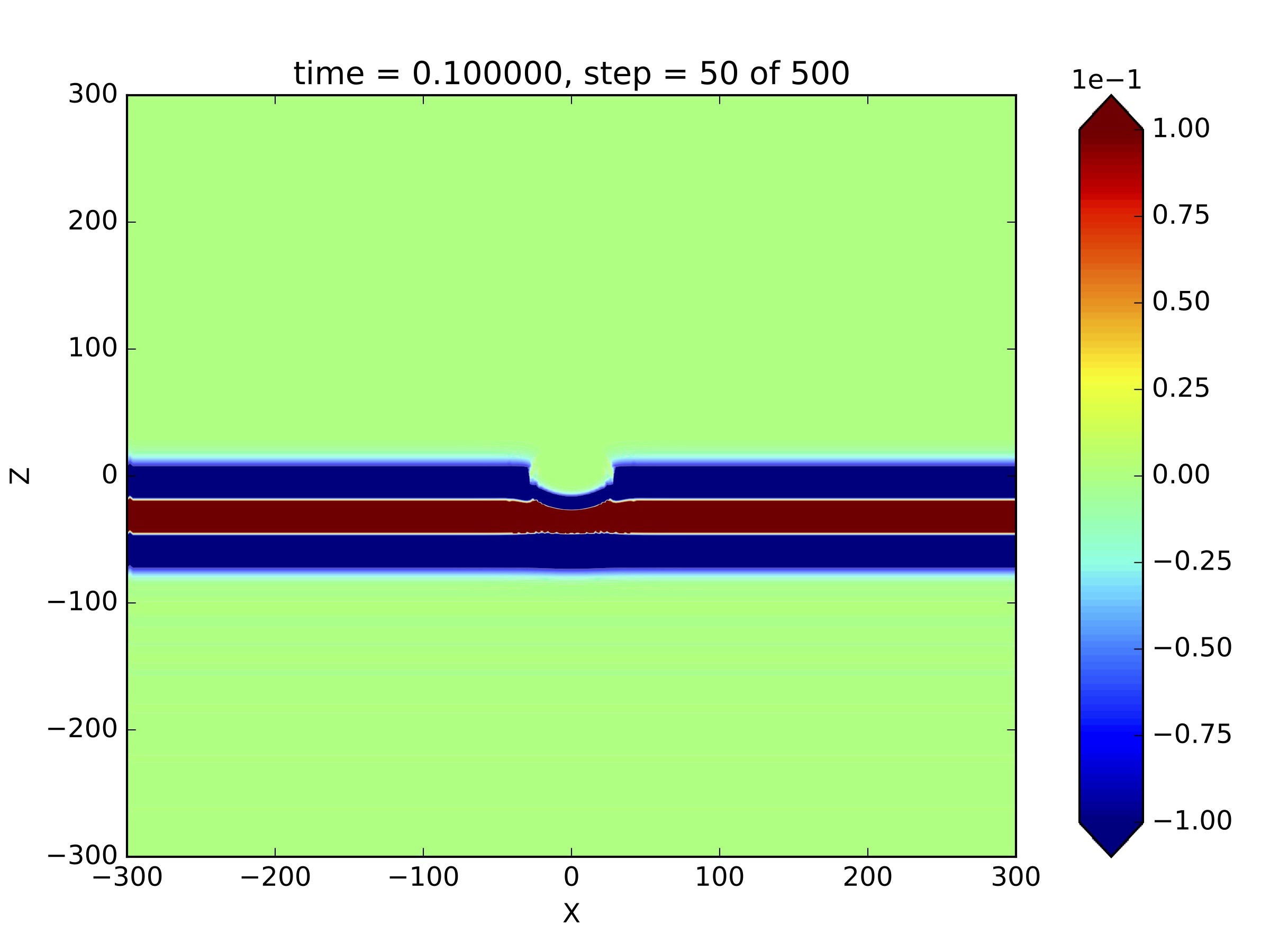}
\includegraphics[width=0.4\textwidth,scale=0.1]{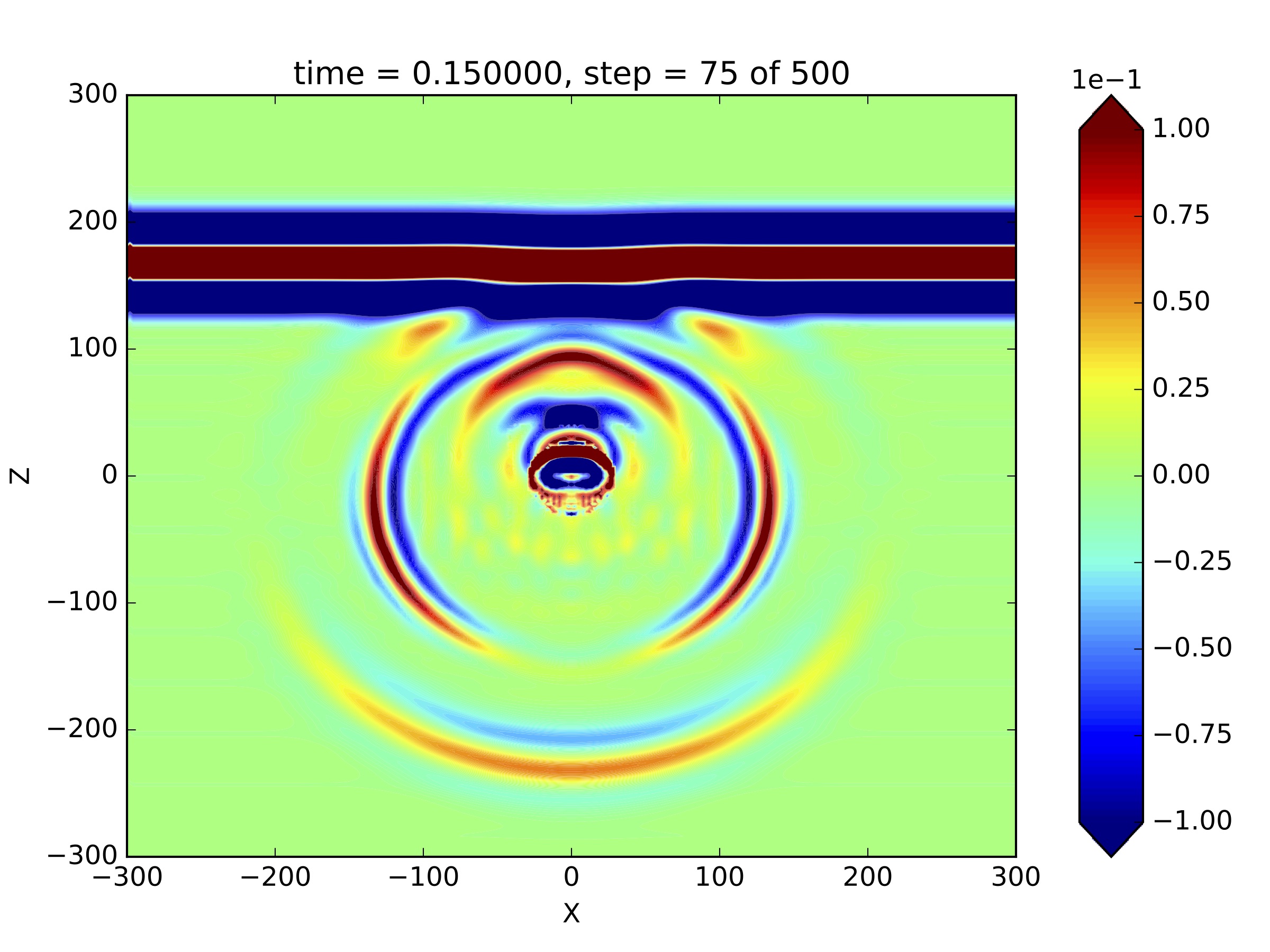}
\includegraphics[width=0.4\textwidth,scale=0.1]{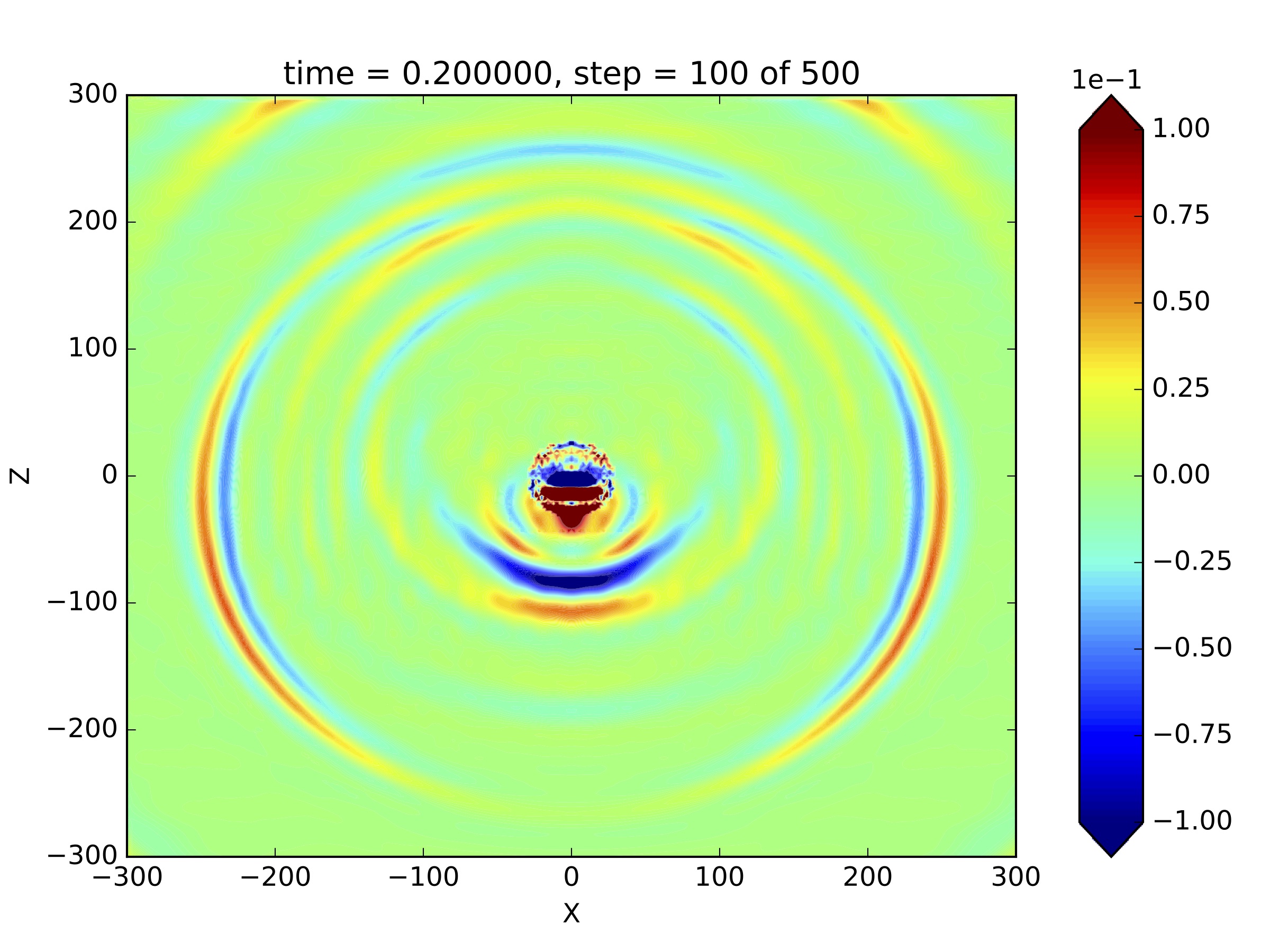}
\includegraphics[width=0.4\textwidth,scale=0.1]{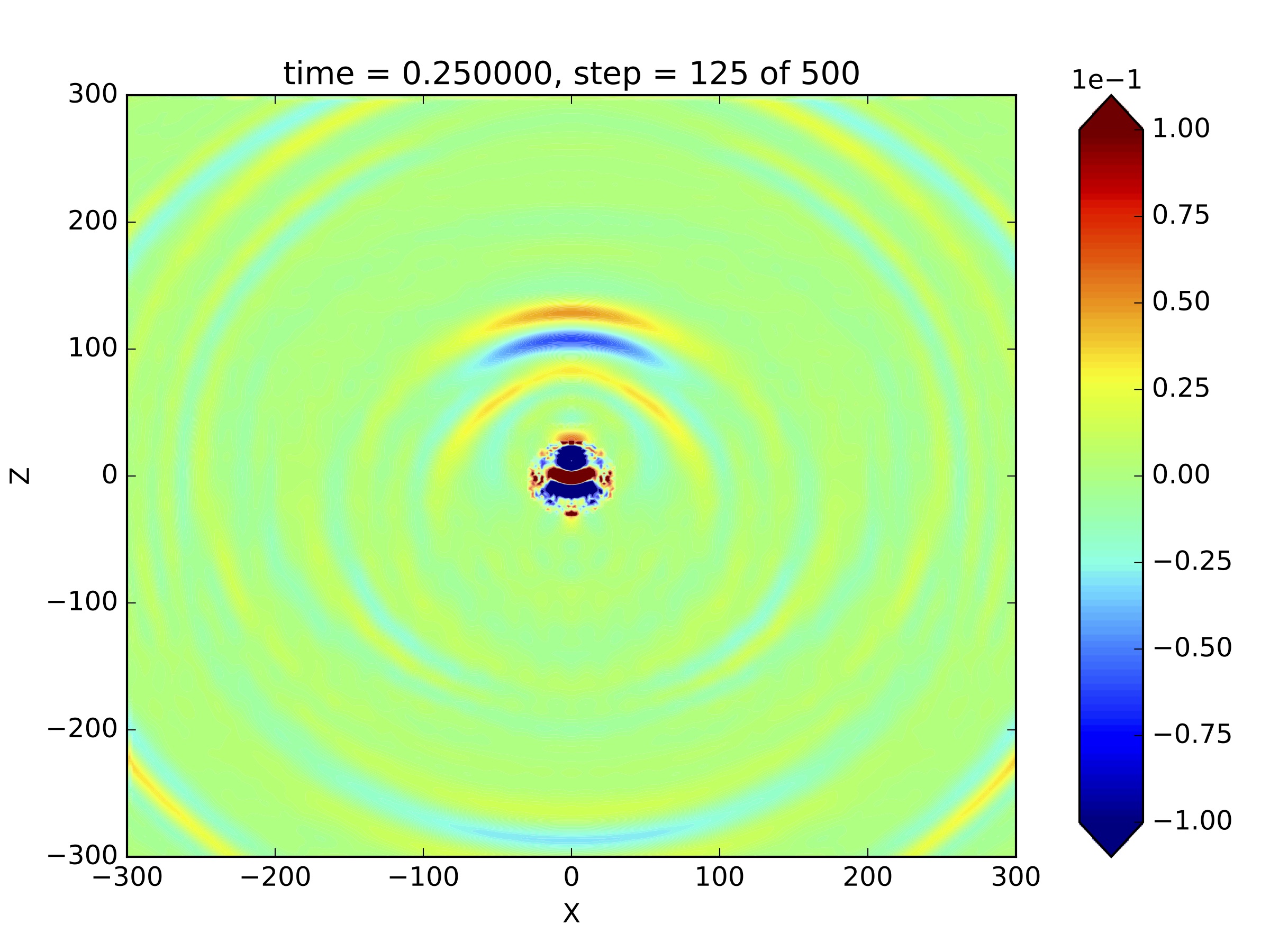}
\includegraphics[width=0.4\textwidth,scale=0.1]{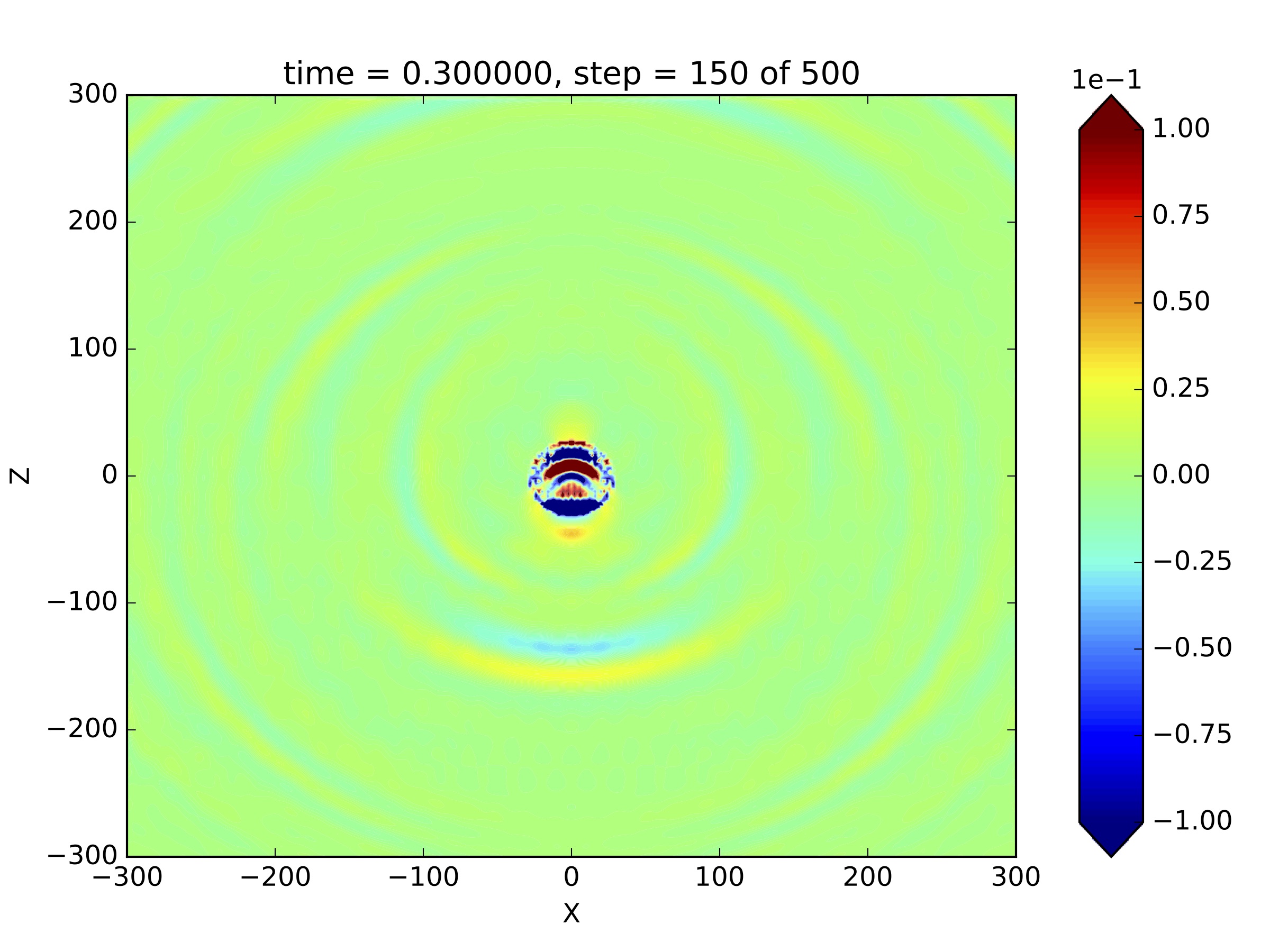}
\includegraphics[width=0.4\textwidth,scale=0.1]{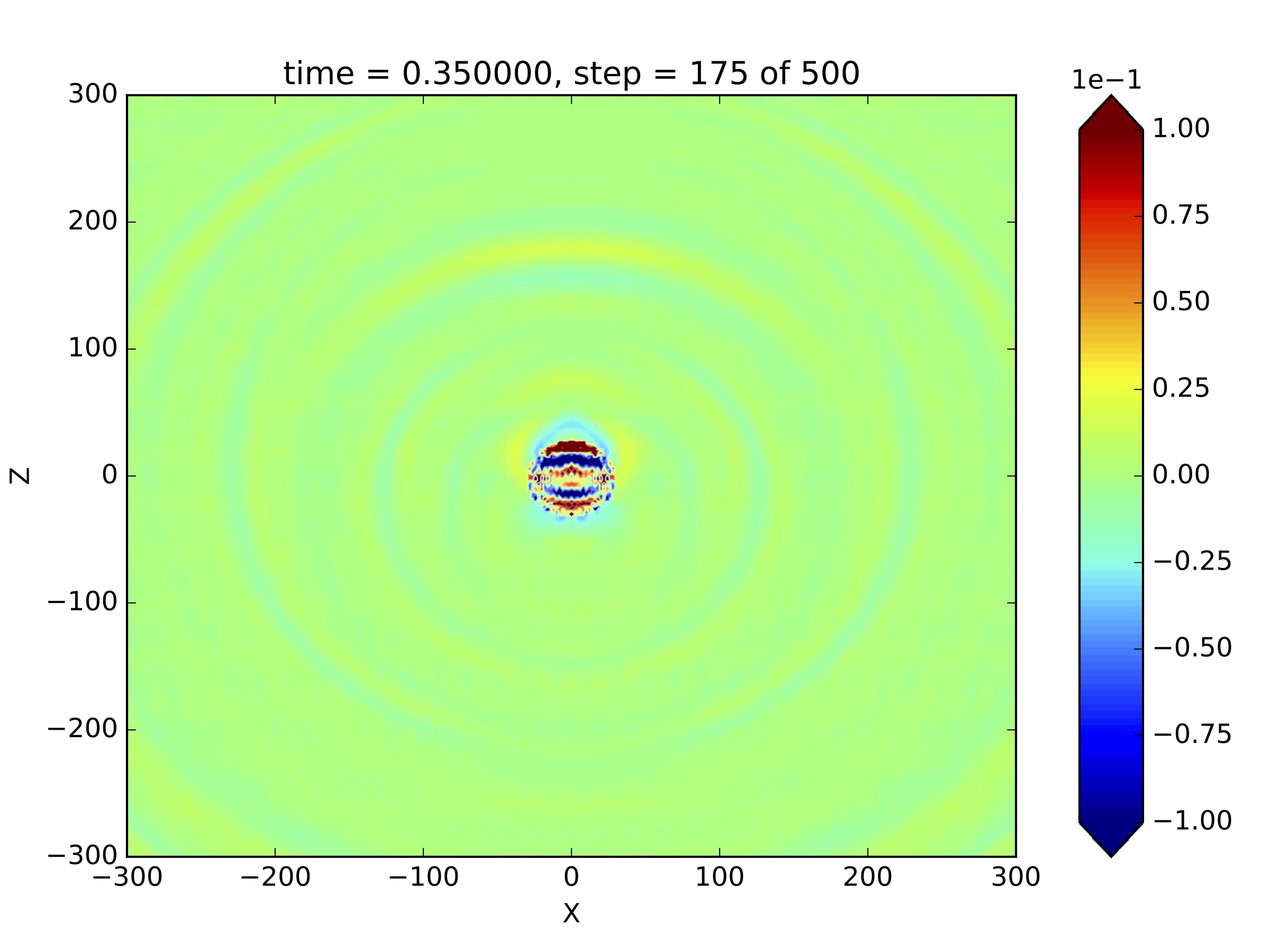}
\includegraphics[width=0.4\textwidth,scale=0.1]{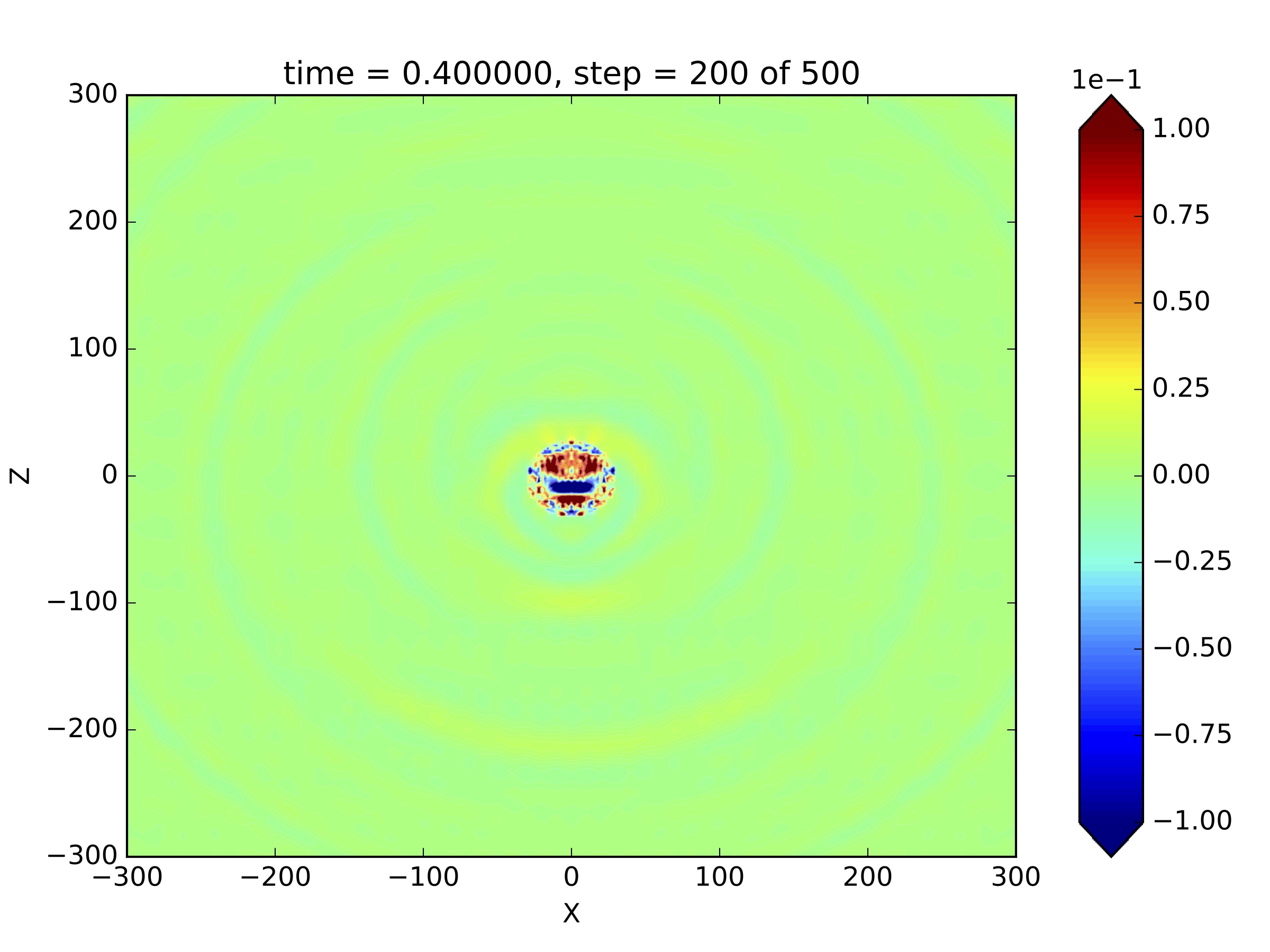}
\hfill
\includegraphics[width=0.4\textwidth]{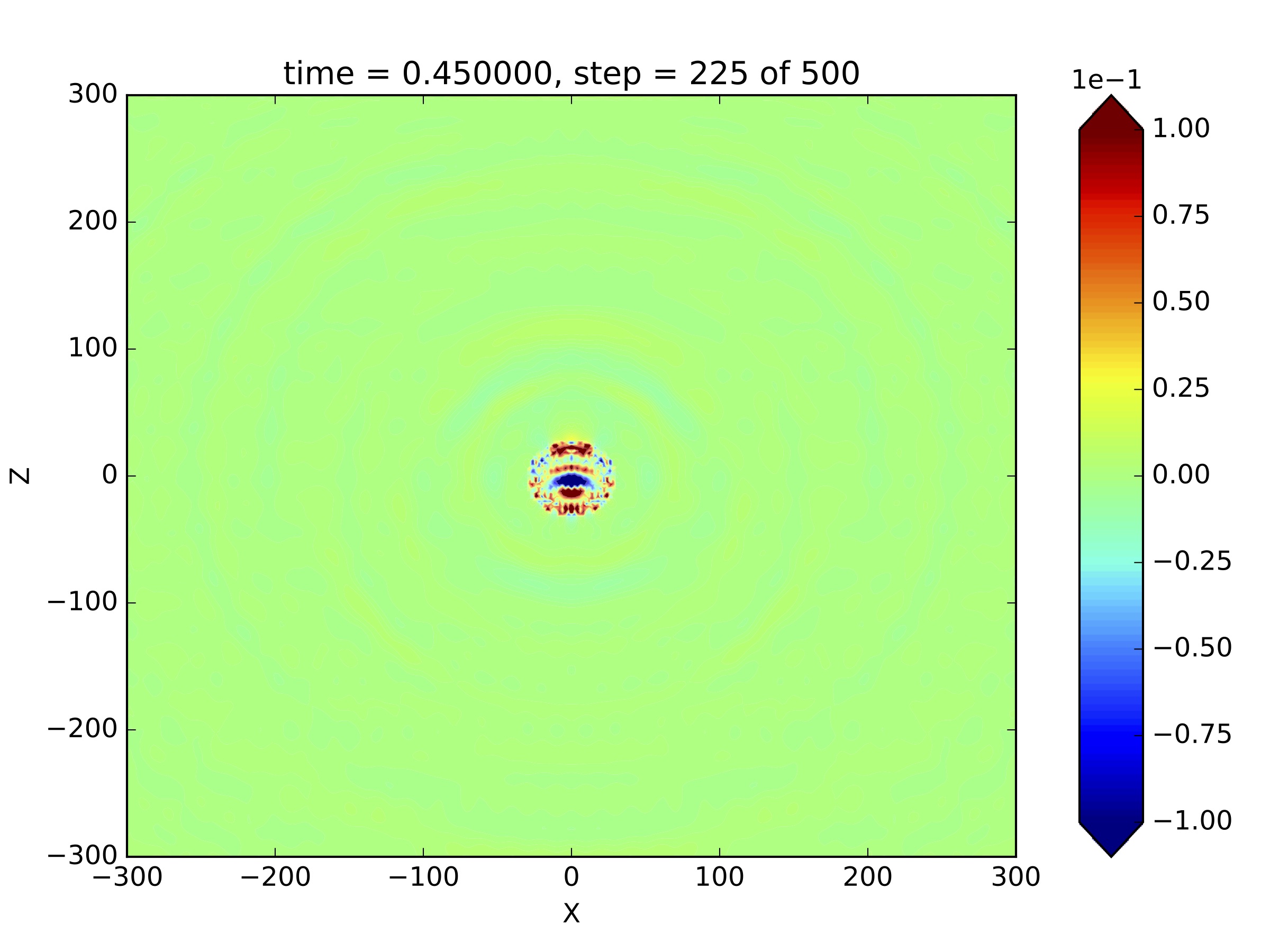}
\caption{Snapshots of the computed wave field.}
\label{fig:snap}
\end{figure}

Profile A in Figure~\ref{fig:profileA} shows the comparison of the displacement time histories in the Z-component (seismograms) for locations crossing the cavity
in vertical direction.  
Seismograms with a gray background show monitoring points inside the acoustic cavity.
The incident plane P-wave reaches the acoustic-elastic interface at about $0.1~s$.
Along this profile only the primary back scattered spherical P- wave (decreasing in amplitude with
distance and time) is visible.  A scattered S-wave is not formed since at $x=0$ the plane
P-wave hits the spherical interface with an incidence angle equal to zero.  
Further, one can see the wave continue to propagate inside the cavity (gray zone) with a lower velocity. Due to the velocity contrast the wave is trapped inside the cavity emitting about 25\% of its energy into the elastic medium each time the wave hits the boundary of the cavity resulting in multiple reverberations decreasing in amplitude with distance and time.
The wave inside the cavity gets also more and more diffracted with time due to the spherical geometry of the cavity. Further, one can see that the incident wave field is
shielded by the cavity creating a shadow zone which causes the suppression of
the incident wave field up to about $10~m$ above the cavity.  Due to 
wave-front healing the incident wave field is present above the cavity and
seems to be unperturbed for $z\gtrsim50~m$. 

The profile A crosses the cavity at the top and bottom, where the direction
normal to the interface points into the $z-$direction.  Thus the $z-$component at
$z=\pm30~m$ in Figure~\ref{fig:profileA} is the normal component with respect
to the interface.  The acoustic-elastic interface condition requires the
continuity of the normal component at the interface which can be seen in
Figure~\ref{fig:profileA}. 
The overall waveform fit is quite satisfying. The direct and multiple
scattered phases are reliably captured.  However some misfit due to numerical
dispersion occurs, getting more pronounced with time.  As discussed in 
section~\ref{sec:input}, this could be overcome by using a finer grid using,
but with a huge computing time.

Figure~\ref{fig:profileB} shows the profile crossing the cavity in horizontal
direction.  Here the seismograms for \mbox{$|x|<30~m$} and \mbox{$|x|\geq30~m$}
are computed in the acoustic and elastic domains, respectively. The incident
wave field passes the profile a little later than $t=0.1~s$. In the elastic
domain the incident field is followed by the primary scattered S-wave. The
transmitted P-wave inside the cavity shows the strongest amplitudes at the
center of the cavity. As described above multiple reverberations occur with decaying amplitudes, cf. Figure 5.

Profile B is oriented horizontally, hence the $z-$components show the tangential
component at the interface.  The physical interface condition only demands
the continuity of the displacement in normal direction.
As discussed in Section~\ref{sec:inter}, the DG-implementation forces all three
components to be continuous at the interface.  However, we can see in
Figure~\ref{fig:profileB} that the discontinuity of the tangential component is
well fitted for the primary transmitted wave.  The scattered numerical waves in the
elastic domain coincide very well with the analytic solution.  The  multiple
internal reverberations are well captured except for numerical dispersion.  A
misfit inside the cavity is present near the interface, which do not seem to
affect the seismograms in the elastic domain.

Figure~\ref{fig:profileCD} shows seismograms for profiles C and D in the back-
and forward-scattered regimes, respectively. The profiles are located in $300~m$
distance above and below the cavity in the elastic medium. In the back-scattered
 regime in Figure~\ref{fig:profileCD} (top) the incident 
plane wave is separated from the scattered waves and arrives earlier in time at
about $t=0.025~s$. The primary scattered P-wave is a distinct wave arrival at
about $t=0.2~s$. Two further scattered P-wave arrivals from internal
reverberations inside the cavity are well pronounced. The first also coincides
with the arrival of the secondary S-wave.  Further, from the cavity decoupled
S-waves occur from internal acoustic reverberations that are from P-to-S
converted during the transmission.  On the $z-$component the spherical-like
scattered P-waves show strong amplitudes near $x=0~m$, while S-waves are more
pronounced for large $|x|$.  At $x=0~m$ S-waves fade out for two reasons, first
no P-to-S conversion takes place for  an incidence angle equal to zero, neither
during reflection nor during transmission of later acoustic reverberations, and
second due to the projection of the shear particle motion on the $z-$direction.
In the forward-scattered regime in Figure~\ref{fig:profileCD} (bottom) the
primary scattered waves directly follow the incident wave. Internal
reverberations cause the later arrivals as discussed above.

All physical features are captured by the numerical solution. The waveform fit
is very satisfying. Small deviations due to numerical dispersion can easily be
overcome by a finer grid or higher polynomial degrees.

\begin{figure}[h!]
\centering
\includegraphics[width=\textwidth, trim= 0cm 0cm 0cm 1.5cm, clip=true]{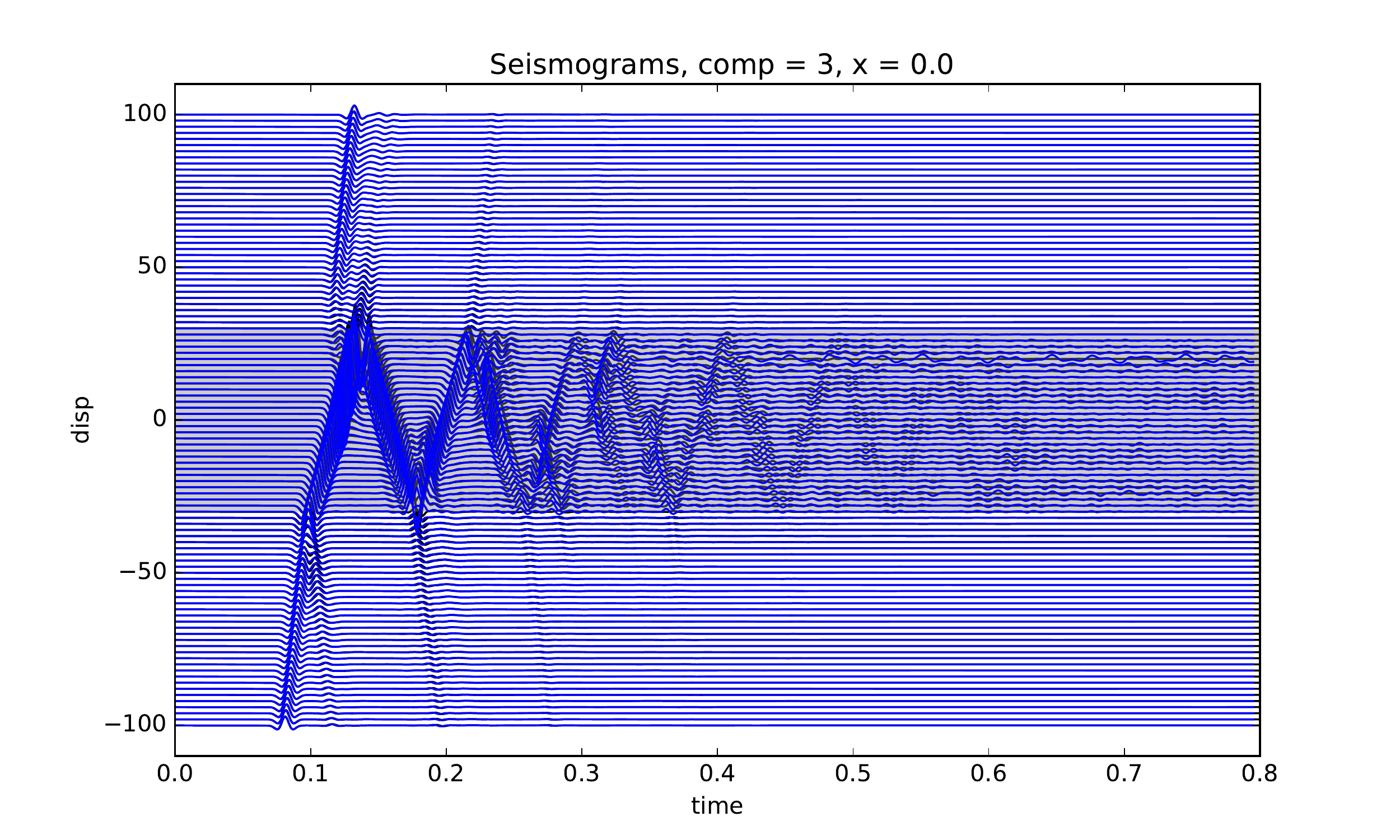}
\caption{Comparison between analytic (black) and numerical (blue) solution along profile A) in Figure~\ref{fig:valid} }
\label{fig:profileA}
\end{figure}

\begin{figure}[h!]
\centering
\includegraphics[width=\textwidth, trim= 0cm 0cm 0cm 1.5cm, clip=true]{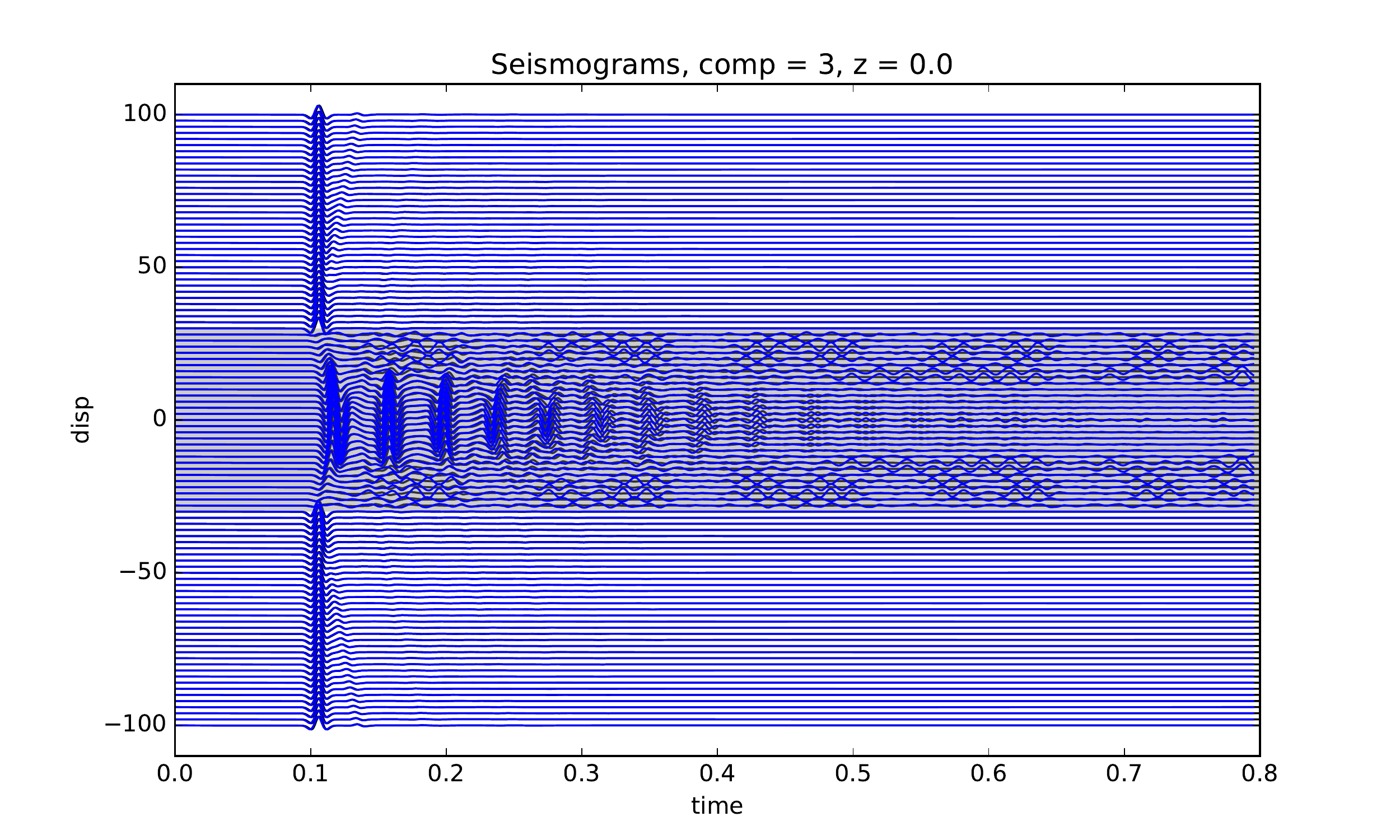}
\caption{Comparison between analytic (black) and numerical (blue) solution along profile B) in Figure~\ref{fig:valid} }
\label{fig:profileB}
\end{figure}

\begin{figure}[h!]
\centering
\includegraphics[width=\textwidth, trim= 0cm 0cm 0cm 0cm, clip=true]{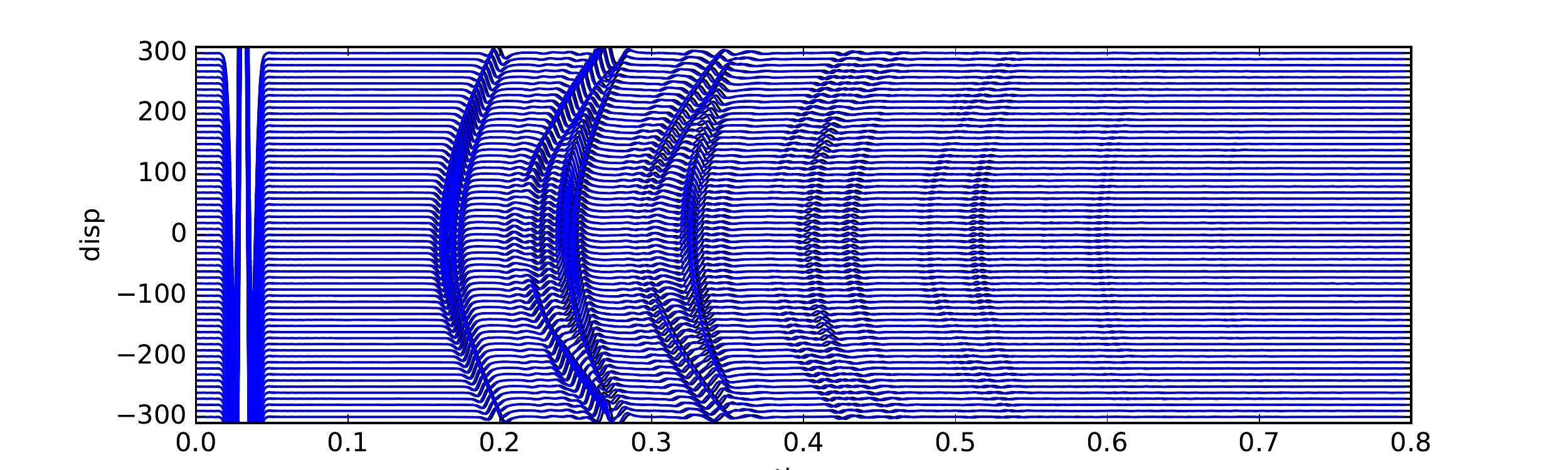}\\
\includegraphics[width=\textwidth, trim= 0cm 0cm 0cm 0cm, clip=true]{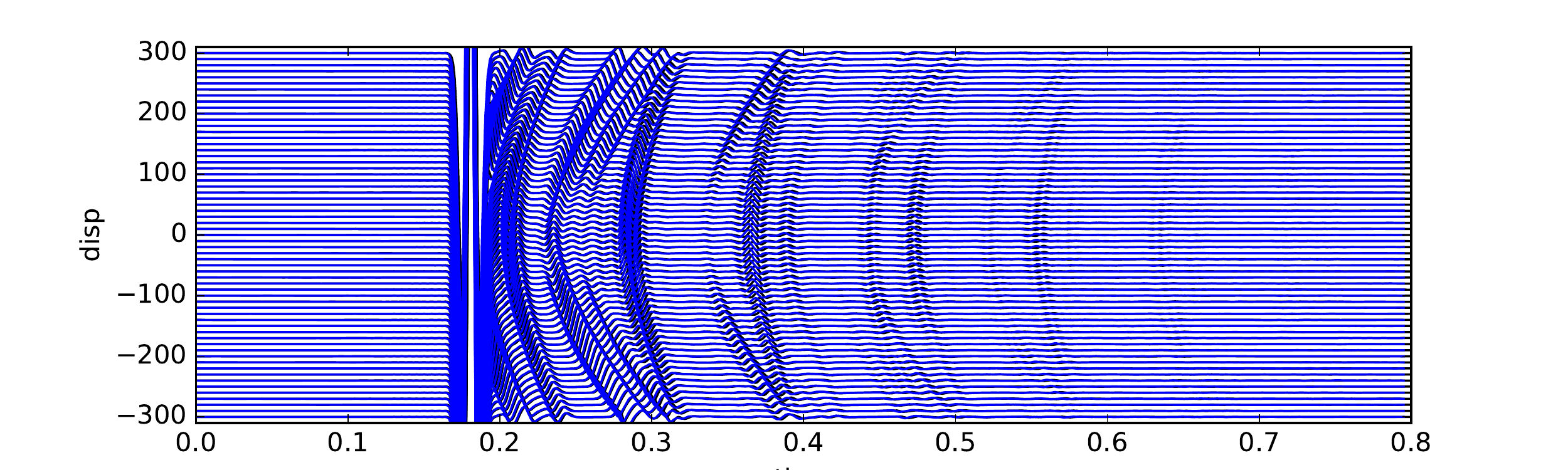}
\caption{Comparison between analytic (black) and numerical (blue) solution along profile C (top) and D (bottom) in Figure~\ref{fig:valid} }
\label{fig:profileCD}
\end{figure}

\FloatBarrier
\section{Conclusion/Discussion}

The numerical computation of the elastic scattering wave field with SPEED stands in good agreement with the analytic solution for a spherical shaped inclusion. The validation for this study paves the way for a comprehensive understanding of the physical characteristic based on numerical computations. 
However, several aspects are open for improvement. { While numerical deviations can be observed in the acoustic domains, all physical features in the elastic domain are well resolved. In the elastic domain no signal information is lost nor any artificial signals are present. }
The acoustic-elastic coupling is an open research topic for implementation in SPEED \cite{Mauri2015} and is currently under investigation. Other combinations of spatial and temporal discretization methods could use local \cite{madec2009energy} or Lax-Wendroff time-stepping \cite{DeBasabe2010Stability}.  An alternative for absorbing boundary conditions was discussed in \cite{Xie2016pml}.

\section*{Acknowledgments}
We acknowledge financial support by the Vienna Science and Technology Fund (WWTF) project MA14-006.
The computational results presented have been achieved in part using the Vienna Scientific Cluster (VSC). Ilario Mazzieri has been partially supported by the research grant Scientific Independence of young Researchers (SIR; RBSI14VT0S); 'PolyPDEs: Non-conforming polyhedral finite element methods for the approximation of partial
differential equations' by Italian Ministry of Education, Universities and Research (MIUR).

\bibliographystyle{abbrv}

\appendix

\section{Comparison of the x-component}

For the sake of completeness we also show here the seismic arrays of the $x$-components along Profile A and B in Figure~\ref{fig:valid}.

Along the vertical profile A the seismic traces simply show no contribution from any shear waves.
Along the horizontal profile B we can again see the multiple reflections inside the cavity and the periodic signals coupling out of the cavity into the surrounding medium.
 
\begin{figure}[h!]
\centering
\includegraphics[width=0.8\textwidth, trim= 0cm 0cm 0cm 0cm, clip=true]{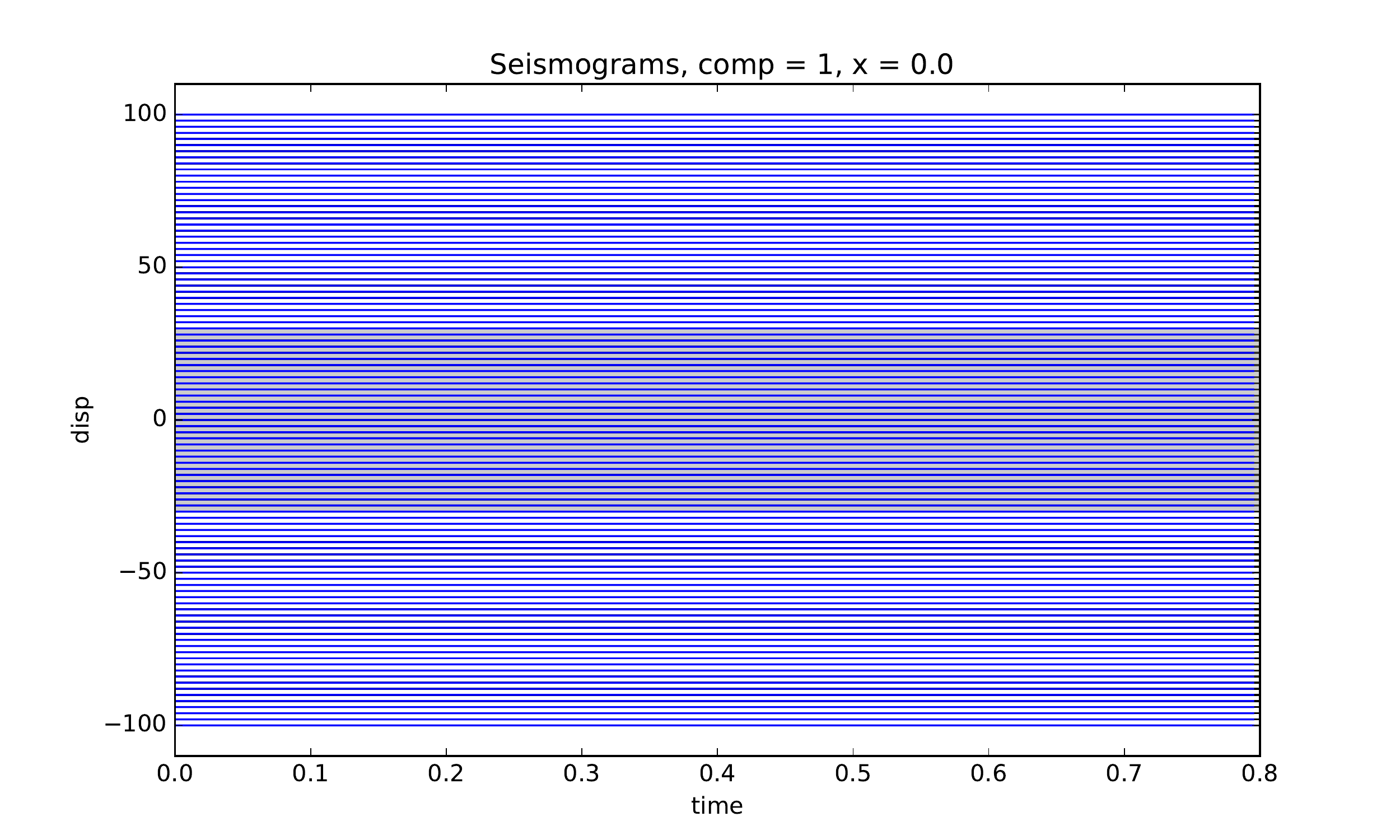}
\caption{Comparison between analytic (black) and numerical (blue) solution along profile A in Figure~\ref{fig:valid}.}
\end{figure}
\begin{figure}[h!]
\centering
\includegraphics[width=0.8\textwidth, trim= 0cm 0cm 0cm 0cm, clip=true]{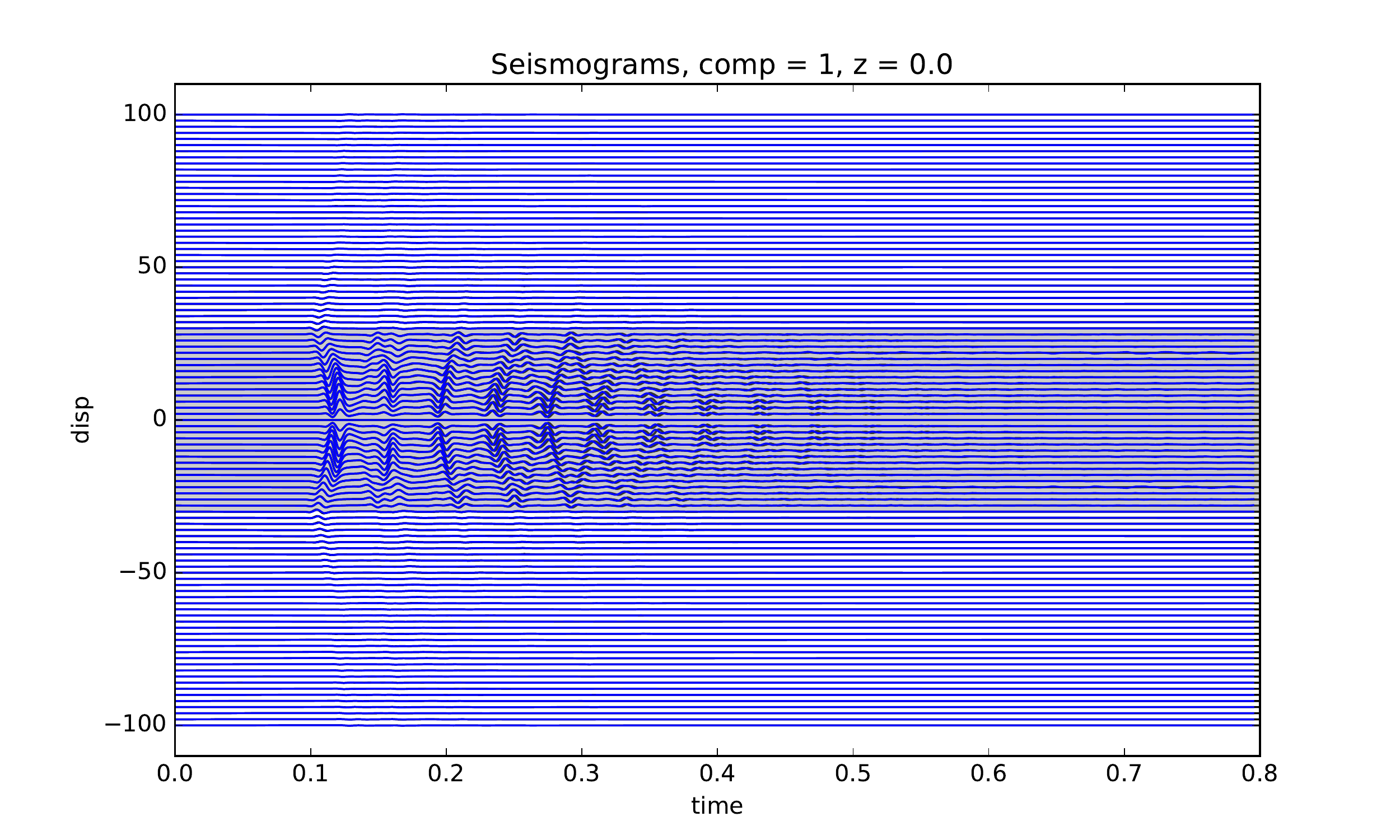}
\caption{Comparison between analytic (black) and numerical (blue) solution along profile B in Figure~\ref{fig:valid}.}
\end{figure}

\end{document}